
\documentstyle{amsppt}
\magnification1200
\NoBlackBoxes

\def\RE{\text{\rm Re }}
\def\phi{\varphi}
\def\L{\fracwithdelims()}

\pageheight{9 true in}
\pagewidth{6.5 true in}

\topmatter
\title Large Character sums
\endtitle
\author Andrew Granville and K. Soundararajan
\endauthor
\author
Andrew Granville and K. Soundararajan
\endauthor
\address{Department of Mathematics, University of Georgia, Athens, 
Georgia 30602, USA} \endaddress
\email{andrew{\@}math.uga.edu} \endemail
\address{Department of Mathematics, Princeton University, Princeton, 
New Jersey 08544, USA} \endaddress
\email{skannan{\@}math.princeton.edu} \endemail
\thanks{The first author is a Presidential Faculty Fellow.  He is also
supported, in part, by the National Science Foundation.  The second 
author is partially supported by the American Institute of Mathematics (AIM).}
\endthanks
\dedicatory  Dedicated to John Friedlander on the occasion of his
45$^{th}$ birthday
\enddedicatory
\endtopmatter

\def\L{\fracwithdelims()}
\def\cbar{\overline{\chi}}

\head {Introduction} \endhead

\noindent  A central problem in analytic number theory is to gain an understanding of  character sums
$$
\sum_{n\le x} \chi(n) ,
$$
where $\chi$ is a non-principal Dirichlet character $\chi \pmod q$.
It is easy to show that such characters sums are always $\leq q$ in
absolute value, while G.~P{\'o}lya and 
I.M.~Vinogradov [3] improved this to $\leq \sqrt{q} \log q$ around 1919,
and H.L.~Montgomery and R.C.~Vaughan [13] to $\ll \sqrt{q} \log\log q$ assuming the
Generalized Riemann Hypothesis (GRH), in 1977.   Up to the constant
this is ``best possible'' since R.E.A.C.~Paley [14] had shown, in 1932,
that there exist characters sums (with real, quadratic characters), that are
$\gg \sqrt{q} \log\log q$.

In many applications one is interested in when the above character sum
is $o(x)$ with $x$ substantially smaller than $q^{\frac 12 +o(1)}$, that is 
$$
\biggl |\sum_{n\le x} \chi(n)\biggr| = o(x). \tag{1}
$$
In 1957, Burgess [2] used ingenious combinatorial methods together with
the ``Riemann Hypothesis for hyperelliptic curves'' to establish (1) 
whenever $x> q^{\frac 14 +o(1)}$,
for any quadratic character mod $q$, with $q$ prime  (and subsequently generalized this to any  non-principal character $\chi \pmod q$ when $q$ is cubefree; with the smaller range $x >q^{\frac 38+o(1)}$ otherwise).
Recently Friedlander and Iwaniec [4] have supplied a different proof of
Burgess's result, and Hildebrand [9] observed that one can 
``extrapolate'' Burgess's bound to the range $x > q^{\frac 14-o(1)}$.  
However, Burgess's range has not been substantially improved 
over the last forty years although it is widely 
believed that such an estimate should hold for $x \gg_\epsilon q^\epsilon$.

In this paper we investigate the distribution of the size of character
sums, and in particular in what range the estimate (1) should hold.  
For example on this question we prove:

\proclaim{Corollary A} Assuming the Generalized Riemann Hypothesis, the
 estimate (1) holds if $\log x/\log\log q \to \infty$ as $q\to \infty$.
This is ``best possible'' in the sense that, for any given $A>0$,
for every prime $q$ there exists a non-principal character $\chi \pmod q$
such $|\sum_{n\leq x} \chi(n)|\gg_A x$ where $x=\log^Aq$.
\endproclaim

  The proof of the first part of this result is inspired by Montgomery and
Vaughan's paper mentioned above. In fact, modifying and refining their
argument 
we will get upper bounds on character sums in all ranges, assuming GRH,
which we believe are close to the truth --- we will discuss a more refined
conjecture below.

  To believe one's upper bounds are close to the truth, one wants to show
that there are character sums of comparable magnitude. Previous arguments
to show that such sums exist, as in Paley's work described above, have
relied in part on using the law of quadratic reciprocity and Dirichlet's
theorem for primes in arithmetic progression to find discriminants for
which many of the small primes are quadratic residues. Such an argument 
seems unlikely to generalize to characters of high order, and might make
one suspicious that perhaps one can only obtain particularly large
character sums (for instance, $\gg \sqrt{q} \log\log q$) 
when the character is real and quadratic. However this is not so, as we
shall show below with a very different proof, involving high moments
of character sums.

  In the large character sums that we exhibited to prove Corollary A,
we showed that they are large by establishing, for those characters, that
the character sum over ``smooth integers'' is particularly large.
Here ``smooth'' refers to integers with only small prime factors, and we define
$$
\Psi(x,y; f) := \sum\Sb n\le x\\ p|n \implies p\le y \endSb f(n) , 
$$ 
for any arithmetic function $f$. Our work on upper and lower bounds motivates 
our belief that character sums can only be large because of 
extraordinary behaviour of the values of $\chi(p)$ for small primes $p$.
We formalize this as the following conjecture:

\proclaim {Conjecture 1}  There exists a constant $A>0$ such that
for any non-principal character $\chi \pmod q$, and for any
$1\le x\le q$ we have, uniformly,  
$$
\sum_{n\leq x} \chi(n) = \Psi(x,y;\chi)+o(\Psi(x,y;\chi_0)),
$$
where $y= (\log q +\log^2 x) (\log \log q)^A$.
\endproclaim

The function $\Psi(x,y):=\Psi(x,y;1)$, the well-known counting function for
smooth numbers, has been extensively investigated. For any  fixed $u>0$,
we know that $\lim_{x\to \infty} \Psi(x,x^{1/u})/x$ exists, and equals 
$\rho(u)$, where $\rho(u)=1$ for $0\leq u\leq 1$ and is 
the real continous function satisfying the differential-delay equation
$\rho^\prime(u)=-\rho(u-1)/u$ for all $u>1$.  We note that 
$\rho(u)=1/u^{u+o(u)}$ as $u\to \infty$.  In \S 3b we will discuss several 
further estimates for $\Psi(x,y)$, though see [10] for a survey. 

Note that Conjecture 1 implies the results of Corollary A, 
and, in fact, further
that if
$$
\Delta(x,q) := \max_{\chi \neq \chi_0} \biggl|\sum_{n\le x} \chi(n)\biggr|
$$
then $\Delta(x,q) \sim \Psi(x,\log q)$ 
whenever $\log x=o( (\log\log q/\log\log\log q)^2)$, for any prime $q$.

Assuming the GRH it is known that there exists 
$n\le \log^{2+o(1)} q$ with $(n,q)=1$ for which $\chi(n)\neq 1$;
assuming Conjecture 1 this would be improved to $n\le \log^{1+o(1)}q$
(see [5] for the latest unconditional work on this problem).

In the wider range $x\le \exp(\sqrt{\log q})$, Conjecture 1 implies that
$$
\biggl| \sum_{n\le x} \chi(n)\biggr| \le \Psi(x,(\log q)^{1+o(1)}) 
= \frac{x}{u^{u+o(u)}}, \qquad \text{where } \qquad u=\frac{\log x}
{\log \log q}. \tag{2}
$$
for any non-principal character $\chi \pmod q$.

We shall establish lower bounds on character sums by various different 
methods in this paper (and in [6] and [7]).  These will imply that,
in most ranges of $x$, the value of  $y$ needs to be at least roughly as large as the value for $y$ given   in Conjecture 1.  

We shall establish that Conjecture 1 holds with 
$y= \log^2 q \log^2 x(\log \log q)^{O(1)}$, assuming GRH, by
extending the method of [13]. This implies the upper bound
$\ll x/u^{u/4+o(u)}$ in (2), as well as the first part of Corollary A.

We shall also establish that Conjecture 1 holds for ``almost all'' characters
$\chi \pmod q$ when $x\le \exp((\log \log q)^{O(1)})$. More
generally we shall show that Conjecture 1 with $y=\log q \log x (\log \log q)^{O(1)}$
holds for almost all non-principal characters $\chi \pmod q$.

Rather than the size distribution, one might be interested in the ``angle
distribution'' of large character sums $\pmod q$. For example, if
a character sum is ``large'', in what directions can it point?  Below
we show, unconditionally, that for any fixed $A>0$, for any given 
angle $\theta$, there are non-principal characters $\chi$ modulo any
prime $q$ for which the character sum up to $\log^Aq$ equals
$\{ e^{i\theta}+o(1)\} \rho(A) \log^A q$. In [7] we show the complementary
result that there are non-principal characters $\chi$ modulo any
prime $q$ for which the character sum up to $q/2$ equals
$\{ e^{i\theta}+o(1)\} (e^\gamma/\pi) \sqrt{q}\log\log q$.

We shall also consider analogues of our results
 for real characters, when appropriate; that is,
$$
\Delta_{{\Bbb R}}(x,q) = \max\Sb q\le |D|\le 2q \\ 
 \endSb \biggl| \sum_{n\le x} \L{D}{n}\biggr|,
$$
where $D$ runs over fundamental discriminants.  
We establish similar and, in some cases, stronger versions of 
the results for $\Delta(x,q)$.

In the next section we give a more technical description of
 our results. In particular our results mostly apply to characters
modulo any integer $q$, not just primes, and with various complicated
error terms.

\head {1. Statement of results} \endhead

We begin with a unconditional, weak version 
of Conjecture 1 which works for  ``almost all'' characters $\chi\pmod q$.

\proclaim{Theorem 1}  Let $1\le x\le q$ be given.  For all but at most 
$q^{1-\frac{1}{\log x}}$ characters $\chi \pmod q$ we have 
$$
\sum_{n\le x} \chi(n) = \Psi(x,y;\chi) + O\biggl(\frac{\Psi(x,y)}{(\log \log q)^2}
\biggr), \qquad \text{whenever } \qquad y\ge\log q\log x (\log \log q)^5.
$$
For all but at most $q^{1-\frac{1}{(\log \log q)^2}}$ characters $\chi\pmod q$
we have 
$$
\biggl|\sum_{n\le x} \chi(n) \biggr| \le \Psi(x, (\log q+\log^2 x)
(\log \log q)^5).
$$
\endproclaim

\remark{Remark} Let $\omega(q)$ denote the number of prime factors
of $q$. Tenenbaum [16] showed that $\Psi(x,y,\chi_0) \asymp (\phi(q)/q) \Psi(x,y)$ whenever $\log y \gg (\log 2\omega(q))(\log\log x)$.
Since $1\geq \phi(q)/q \gg 1/\log\log q$ we see that the 
error term in Theorem 1 can be rewritten as 
$O( \Psi(x,y,\chi_0)/\log\log q)$ in this range.
\endremark

Assuming the GRH we can establish results similar to (but weaker than) Theorem 
1, but valid for all non-principal characters.  The prototype for 
our result appears as Lemma 2 in [13].  There, Montgomery and Vaughan 
show that if $\chi \pmod q$ is non-principal and the GRH holds then
$$
\sum_{n\le x} \chi(n) = \Psi(x,y) + O(xy^{-\frac 12} \log^4 q),
$$
when $\log^4 q \le y\le x\le q$.  Their objective was not to establish this
in as wide a range as possible; however, ours is, so 
we modify and refine their method to 
obtain the following result.

\proclaim {Theorem 2} Assume that the Generalized Riemann Hypothesis
holds true, and let $\chi$ be any non-principal character $\pmod q$.
If $1\le x\le q$ and $y\geq \log^2q\log^{2} x (\log \log q)^{12}$  then
$$
\sum_{n\leq x} \chi(n) = \Psi(x,y;\chi)+O\biggl( \frac{\Psi(x,y)}
{(\log \log q)^2}\biggr).
$$
Further
$$
\biggl| \sum_{n\leq x} \chi(n) \biggr| \ll
\Psi(x, \log^{2}q (\log \log q)^{20}),
$$
and so the estimate (1) holds when $\log x/\log\log q \to \infty$ as $q\to \infty$.
\endproclaim
\remark{Remarks} To compare this with Montgomery and Vaughan's result, the 
error term in the first part of Theorem 2 could have been written as the
rather more complicated
$O( \Psi(x,y) \log q \log x (\log\log q)^4/\sqrt{y})$. Similarly,
the error term in the first part of Theorem 1 can be considerably 
sharpened.

As in Theorem 1 the error term can be rewritten as 
$O( \Psi(x,y,\chi_0)/\log\log q)$ when
 $\log y \gg (\log 2\omega(q))(\log\log x)$.
\endremark

We now proceed to the problem of finding large character sums, beginning
with the range $x\le \exp((\log \log q)^{2-\epsilon})$.  Here we 
get large character sums, pointing in any given direction.

\proclaim{Theorem 3} Let $q$ be large, and suppose 
$\log x\le \frac{(\log \log q)^2}{(\log \log \log q)^2}$.  
For all $|\theta|\le \pi$ 
there are at least $q^{1-\frac{2}{\log x}}$ characters 
$\chi \pmod q$ for which 
$$
\sum_{n\leq x} \chi(n) =   e^{i\theta}\ \Psi(x,\log q;\chi_0)  
+ O\biggl( \Psi(x,\log q) 
\biggl( \frac{1}{\log x}+\frac{\log x (\log \log \log q)^2}{(\log \log q)^2}\biggr)\biggr).
$$
If $q$ has no prime factors below $\log q$ then we may write 
the above as 
$$
\sum_{n\leq x} \chi(n) =   x e^{i\theta}\ \rho\biggl(\frac{\log x}{\log \log q}
\biggr)  \biggl(1+
O\biggl( \frac{1}{\log x}+\frac{\log x (\log \log \log q)^2}{(\log \log q)^2}\biggr)\biggr).
$$
\endproclaim

This implies the second part of Corollary A.

Theorem 3 is not useful when $q$ has many prime factors below $\log q$.
We next deduce, by a very different method,
 lower bounds of more or less the same strength for these cases.

\proclaim{Theorem 4}  Suppose $x=(10\log q)^{B} = q^{o(1)}$ for some $B\ge 1$. 
Then 
$$
\max_{\chi \neq \chi_0} \biggl| \sum_{n\le x} \chi(n) \biggr| 
\gg \frac{x^{\frac{1}{2}+ \frac{[B]}{2B}}}{(4\log x)^{[B]}}.
$$
If, in addition, $q$ has less than $(\log q)^{B/(B+1) -\epsilon}$ distinct 
prime factors then 
$$
\max_{\chi \neq \chi_0} \biggl| \sum_{n\le x} \chi(n) \biggr| 
\gg \frac{x}{(4\log x)^{[B]+1}}.
$$
\endproclaim

Applying Theorem 4 appropriately, we can deduce the following 
corollaries.  

\proclaim{Corollary 1}  If $\log x \ge (\log  \log q)^2$ then 
$$
\max_{\chi \neq \chi_0 } \biggl| \sum_{n\le x} \chi(n) \biggr| 
\gg x \exp\biggl( -(1+o(1)) \log x \frac{\log \log x}{\log \log q}\biggr).
$$
If, in addition, $q$ has less than $(\log q)^{1-\epsilon}$ distinct 
prime factors then this bound holds in the extended range $\log x/\log \log q
\to \infty$.  
\endproclaim

\remark{Remark} There are $\ll q/\exp( (\log q)^{1-\epsilon} )$ integers $q\leq x$ failing the restriction ``$q$ has less than $(\log q)^{1-\epsilon}$ distinct 
prime factors''.
\endremark

\proclaim{Corollary 2}  Fix $\sigma$ in the range $\frac 12\le \sigma<1$. 
If $(\log q)^{\frac 1{1-\sigma}} 
\le x\le \exp((\log q)^{1-\sigma +o(1)})$ then 
$$
\max_{\chi \neq \chi_0} \biggl|\sum_{n\le x} \chi(n) \biggr| 
\gg x^{\sigma}.
$$
If, in addition, $\omega(q) \le (\log q)^{\frac 12-\epsilon}$ 
then this bound holds 
whenever $x\ge (\log q)^{1+\epsilon}$.  In any case we have
$$
\max_{x \ge 1}\max_{\chi \ne \chi_0} \frac{1}{x^\sigma} \
\biggl|\sum_{n\le x} \chi(n)\biggr| \gg 
\exp\left( \frac{(\log q)^{1-\sigma}}{14\log\log q} \right).
$$
\endproclaim

So far we have dealt with the range $x\le \exp((\log q)^{\frac 12-\epsilon})$.
We now proceed to the range when $x$ is larger, dealing first with the 
range $\log \log x = (\frac 12+o(1))\log \log q$.

\proclaim{Theorem 5}  Suppose that $\log x = \tau \sqrt{\log q \log\log q}$
with $\tau = (\log\log q)^{O(1)}$, and let $\eta=\tau+1/\tau$. There
exists a constant $c>0$ such that for any sufficiently large $q$, 
there exists a non-principal character  $\chi\pmod{q}$ for which 
$$
\frac{1}{\sqrt{x}}\biggl| \sum_{n\leq x} \chi (n) \biggr| 
 \gg  \exp \biggl( c\frac{1+\log (\eta\tau)}
{\eta}\sqrt{\frac{\log q}{\log\log q}} 
\biggr) .
$$
\endproclaim

As a consequence we get Corollary 3 below, which improves Corollary 2 
in the case $\sigma =\frac12$.

\proclaim{Corollary 3}  There exists a constant $c > 0$ such that
for all integers $q$
$$
\max_{\chi \ne \chi_0} \max_{x \ge 1} \frac{1}{\sqrt{x}} \ 
\biggl| \operatornamewithlimits\sum\limits_{n\le x}
\chi(n)\biggr| \gg \exp\biggl(c\sqrt{\frac{\log q}{\log\log q}}\biggr) .
$$
\endproclaim

Next we consider 
the range when $\log x/\sqrt{\log q \log \log q}$ is large, but $x$ is
smaller than $q^{\epsilon}$.  

\proclaim{Theorem 6}  Suppose both $\log q/\log x$ and $\log x/\sqrt{\log q
\log\log q} \rightarrow \infty$.  There exists a non-principal character 
$\chi\pmod{q}$ for which 
$$
\frac{1}{\sqrt{x}} \biggl| \sum_{n\leq x} \chi (n) \biggr| 
 \gg \biggl(\frac{\log x}{\sqrt{\log q \log\log q}}\biggr)^{(1 + o(1))
\frac{\log q}{\log x}}.
$$
\endproclaim

When $x$ is as large as a power of $q$ we obtain:

\proclaim{Theorem 7}   Let $k\geq 2$ be an integer and suppose 
$\exp(\frac{\log q}{\log \log q})\le x<q^{\frac 1k}$.  Then
there exists a non-principal character $\chi\pmod{q}$ for which 
$$
\frac{1}{\sqrt{x}} \biggl| \sum_{n\leq x} \chi (n)\biggr| \gg_k  \ 
(\log q)^{\frac{(k -1)^2}{2k} + o(1)}.
$$
\endproclaim

Once $x\ge q^{\frac 12}$ Theorem 7 reduces to the bound $\Delta(x,q)
\ge \sqrt{x}(\log q)^{o(1)}$ which follows immediately from the 
mean square of $\sum_{n\le x} \chi(n)$.  However it is possible to
obtain non-trivial information here by appealing to (essentially) the Poisson 
summation formula.
We quote P{\' o}lya's Fourier expansion (see Lemma 1 of [13])
$$
\sum_{n\le x} \chi(n) = \frac{\tau(\chi)}{2\pi i} \sum\Sb h=-H \\ 
h\neq 0\endSb^{H} \frac{\overline{\chi}(h)}{h} 
(1-e(-\tfrac{hx}{q})) + O(1+qH^{-1}\log q),\tag{3}
$$
where $\chi$ is primitive, and $\tau(\chi)$ is the usual Gauss sum.   
Since $|\tau(\overline{\chi})|
=\sqrt{q}$, (3) suggests a relation of the type $\Delta(x,q) $`$=$'$
\frac{x}{\sqrt{q}}\Delta(\frac qx,q)$; now $\frac qx \le q^{\frac 12}$ so 
that applying the ideas behind our earlier Theorems should lead to 
a good lower bound for $\Delta(x,q)$.  While we cannot 
show such a result for every $x$, using (3) we can obtain good 
bounds for $\Delta(t,q)$ for some $t\le x$.  Naturally 
one would expect $\Delta(t,q)$ to be an increasing function of $t$ (at least
most of the time) but we don't know how to prove this.  For convenience, 
we state this result only for primes $q$, so that every non-principal 
character is primitive.  
\proclaim{Theorem 8} Let $q$ be a large prime.  
Given $\exp(c\sqrt{\log q}) \ge N \ge 2$ (for a small positive constant $c$) we have 
$$
\max_{t\le q/N} \ \max_{\chi \neq \chi_0 \pmod q} \ 
\biggl|\sum_{n\le t} \chi(n)\biggr| \gg \sqrt{q} \ \frac{1}{N}
\Psi\biggl( N, \frac{\log q}{(\log \log q)^{10}}\biggr).
$$
When $\log N = \tau \sqrt{\log q \log \log q}$ with $\tau =
(\log \log q)^{O(1)}$ we have (for a small positive constant $c$ and 
$\eta=\tau+1/\tau$)
$$
\max_{t\le q/N} \ \max_{\chi \neq \chi_0 \pmod q} 
\ \biggl| \sum_{n\le t} \chi (n) \biggr| \gg
\sqrt{q/N} \ \ \exp \biggl( c \frac{1+\log (\eta \tau)}{\eta} 
\sqrt{\frac{\log q}{\log \log q}}\biggr).
$$
If both $\log q/\log N$ and $\log N/\sqrt{\log 
q\log \log q} \to \infty$ then 
$$
\max_{t\le q/N} \ \max_{\chi \neq \chi_0 \pmod q} 
\ \biggl| \sum_{n\le t} \chi (n) \biggr| \gg  
\sqrt{q/N} \ \ \biggl(\frac{\log N}
{\sqrt{\log q \log \log q}} \biggr)^{(1+o(1))\frac{\log q}{\log N}}.
$$
Lastly if $\exp(\frac{\log q}{\log \log q}) \le N \le 
q^{\frac{1}{k}-\epsilon}$ for an integer $k\ge 2$ then 
$$
\max_{t\le q/N} \ \max_{\chi \neq \chi_0 \pmod q} 
\ \biggl| \sum_{n\le t} \chi (n) \biggr| \gg_k \sqrt{q/N} \ \ 
(\log q)^{\frac
{(k-1)^2}{2k}+o(1)}.
$$
\endproclaim

Several different authors (for example [1]) gave the same
explicit version of Paley's result: There are infinitely many non-square, positive integers $q$, and 
integers $x=x_q$ for which
$$
\biggl|\sum_{1 \le n \le x} \left(\frac{q}{n}\right)\biggr| \gtrsim
\frac{e^\gamma}{\pi}\sqrt{q} \log\log q, \tag{4}
$$
where $\gamma \approx 0.5772156649\dots$ is the Euler-Mascheroni constant.
We can prove that there are many 
characters $\chi \pmod q$ for which $\sum_{n\le q/2}\chi(n)$ is of 
such large magnitude, and points in any given direction, for any given
prime $q$.  Further, 
whenever $q (\log q)^{-A} \le x\le q$ we can show that $\Delta(x,q) 
\gg  \rho_A \sqrt{q} \log \log q$, where $\rho_A=1/A^{A+o(A)}$ as $A\to \infty$.  The proofs of these results will appear in [7],
because they are more closely related to the methods of that 
paper.  Note, though, in Theorem 11 below we obtain some results
of this type for real characters.

We now turn our attention to getting bounds for $\Delta_{\Bbb R}(x,q)$: that is, exhibiting
large character sums for real characters. We begin by showing that the lower 
bound implicit in Conjecture 1 holds in a very wide range for real characters.

\proclaim {Theorem 9} Suppose $q$ is large, and that $1\le x\le 
\exp(\sqrt{\log q})$.  Then
$$
\max\Sb q\le |D|\le 2q \endSb 
\biggl|\sum_{n\leq x} \L{D}{n}\biggr|
\geq \Psi(x, \tfrac 13\log q) .
$$
Consequently for a fixed real number $B$ there are fundamental discriminants 
$D$ in the range $q\le |D|\le 2q$ with 
$$
\sum_{n\le x} \L{D}{n} \ge (\rho(B) +o(1)) x \gg_B x, \qquad
\text{where     } x=(\tfrac 13\log q)^B
$$
\endproclaim

Theorem 9 is the analogue of Theorems 3 and 4 above.  From 
Theorem 9 we can deduce the analogues of Corollaries 1 and 2 for 
real characters. 

It seems to have been widely believed that ${\mathop{\sum}_{N<n\le N+x}} 
\L{D}{n} = o(x)$
when $x/\log^2D\to \infty$ (see, for instance, page 379 of [11]), 
perhaps in analogy with the known result ${\mathop{\sum}_{p\le x}} \L{D}p  
= o(\pi(x))$ in this range, assuming GRH.  
However Theorem 9 shows that this widely held
view is false. It seems safe to hazard the guess that, for
all non-principal characters $\chi \pmod q$ we have, uniformly,
$$
\sum_{n=N}^{N+x}  \chi(n) \ll x^{1-1/\log\log q} . 
$$

Set 
$$
\alpha(B) = \limsup\Sb |D| \to \infty\endSb 
 \frac{1}{(\log |D|)^B}\sum_{n\le (\log |D|)^B} 
\L{D}{n}. 
$$
Clearly $\alpha(B)=1$ for $0\le B\le 1$, and from Theorem 9, 
we know that $\alpha(B)\ge \rho(B)$. If the GRH is true  
then $\alpha(B) \le \rho(\frac B2)$, by Theorem 2.  Conjecture 1 predicts 
that $\alpha(B)=\rho(B)$ but this is not known for any $B> 1$.  
In Theorem 3 we obtained large character sums pointing in any given direction.
Mark Watkins asked us if the analogue for real character sums holds: 
that is, can one get real character sums to be large and negative? 
Precisely, what can one say about 
$$
\beta(B) := \liminf\Sb |D|\to \infty\endSb 
 \frac{1}{(\log |D|)^B}\sum_{n\le (\log |D|)^B} 
\L{D}{n}.
$$
Interestingly $\beta(B)$ can never be as small as $-1$.  
Indeed in [6] 
we have shown that 
$\beta(B) \ge \delta_1 =-0.656999\ldots$ (see Theorem 1 of [6]
 for a definition of $\delta_1$) for all $B$, and in fact 
$\beta(B) =\delta_1$ for $0\le B\le 1$.  Answering Watkins' 
question we also show there that $\beta(B)<0$ for all $B$, 
but it is an open problem to determine $\beta(B)$ and $\alpha(B)$ for $B>1$.

We obtain the 
following analogue of Theorems 5, 6 and 7, but in a much wider range.

\proclaim{Theorem 10}  Suppose that  $q$ is large, and 
$\exp ( (\log q)^{\frac 12})\le x\le q/\exp ( (\log q)^{\frac 12})$.
Then there exist fundamental discriminants 
$D$ in the range $q\le |D| \le 2q$ with
$$
\frac{1}{\sqrt{x}} \sum_{n\leq x}\L{D}{n} \gg 
\exp \biggl( (1+o(1)) \frac{\sqrt{\log q }}{\log\log q} \biggr) .
$$
\endproclaim

Notice that Theorem 10 is much stronger than the bounds of 
Theorems 6 and 7, as soon as $\log x\ge \sqrt{\log q} (\log \log q)^2$.  
This difference is especially noticeable when $x$ is like a small power of
$q$, and suggests that Theorems 6 and 7 are unlikely to be 
``best possible.'' 

In the next result we use Poisson summation (as discussed after (3))
to get lower bounds for character sums when $x$ is very large, in
terms of smooth numbers. This suggests that we should be able to make
another conjecture like Conjecture 1 for large $x$, which takes this
natural symmetry into account. We have not yet felt able to 
formulate this appropriately.

\proclaim{Theorem 11}  Let $q$ be large. For any 
$\exp(\sqrt{\log q}) \ge N \ge 2$ there exists a fundamental 
discriminant $D$ with $q\le |D|\le 2q$ such that
$$
\biggl|\sum_{n\le |D|/N} \L{D}{n} \biggr| 
\gg \sqrt{q} \ \frac{1}{N} \ \Psi \left( N,  \frac 19 \log q\right)
\left\{ 1 + \frac{\log \log q}{\log (A+2)} \right\}, 
$$
where $N=(\tfrac 19\log q)^A$.
In particular if $\exp((\log \log q)^2)\ge N \ge 2$ then
there exists a fundamental 
discriminant $D$ with $q\le |D|\le 2q$ such that
$$
\biggl| \sum_{n\le |D|/N} \L{D}{n}\biggr|  
\gg \frac{\rho(A)}{\log (A+2)}\ \sqrt{|D|} \log \log |D| .
$$
\endproclaim

\head 2.  The plan of attack \endhead

\noindent We define complex, multiplicative random variables $X_n$ as follows:
$X_n$ is multiplicative: that is, if 
$n=\prod_i p_i^{a_i}$ then $X_n =\prod_{i} X_{p_1}^{a_i}$.  For primes $p$,
$X_p$ is equidistributed on the unit circle, and for different
primes $p$ and $q$, $X_p$ and $X_q$ are independent.  Thus 
${\Bbb E}(X_m \overline{X_n}) = 1$ if $m=n$, and ${\Bbb E}(X_m {\overline X_n})
=0$ otherwise.  Here, and below, ${\Bbb E}(\cdot)$ denotes the expectation.


Let $f$ be any arithmetical function, and $k$ and $n$ be integers.  Below 
we shall put
$$
d_{k,f}(n,x) = \sum\Sb m_1\ldots m_k = n\\ m_i\le x\endSb f(m_1) \ldots
f(m_k),
$$
so that 
$$
\biggl(\sum_{n\le x} \chi(n) f(n)\biggr)^k =\sum_{n\le x^k} d_{k,f}(n,x) 
\chi(n), \qquad \text{and     }\biggl(\sum_{n\le x} X_n f(n)\biggr)^k = 
\sum_{n\le x^k} d_{k,f}(n,x) X_n.
$$
We shall abbreviate $d_{k,f}(n,x)$ to $d_k(n,x)$ when $f$ is the 
function $f(n)=1$.

\proclaim{Lemma 2.1}  Let $x$, $q$, and $k$ be integers with
$x^k \le q$; and let $f$ be any arithmetic function.  Then 
$$
\align
\frac{1}{\phi(q)} \sum_{\chi\pmod q} \biggl| 
\sum_{n\le x} \chi(n)f(n)\biggr|^{2k}
= \sum\Sb n\le x^k\\(n,q)=1\endSb |d_{k,f}(n,x)|^2
&={\Bbb E}\biggl( \biggl| \sum\Sb n\le x\\ (n,q)=1\endSb X_n f(n)\biggr|^{2k}
\biggr).\\
\endalign
$$
\endproclaim
\demo{Proof}  This is immediate from the definition of $X_n$, and 
the orthogonality of the characters $\pmod q$:
$$
\frac{1}{\phi(q)} \sum_{\chi\pmod q } \chi(a) {\overline \chi}(b) 
= \cases
1 & \text{ if } a\equiv b\pmod q, \ \ (ab,q)=1\\
0& \text{  otherwise}.\\
\endcases
$$
\enddemo  

Our plan (see \S 4, and \S 6) is to obtain large lower bounds for the quantity 
in Lemma 2.1 (in the case $f(n)=1$) so as to 
obtain large non-trivial character sums.  In order to do this, we need to 
eliminate the principal character term (which is often large for 
trivial reasons) which is included in the sum in Lemma 2.1.

For any arithmetic function $f$ we define
$$
\Delta_f(x,q) := \max_{\chi \neq \chi_0 \pmod q} 
\biggl| \sum_{n\le x} \chi(n) f(n)\biggr| .
$$

\proclaim{Proposition 2.2} Let $q$ be large, $x\ge \log q$, and 
suppose $k$ is an integer with $x^k \le q$.  For  any 
arithmetic function $f$ we have
$$
\Delta_f(x,q)^{2k} \gg \frac{1}{\phi(q)} \sum_{\chi \pmod q} \biggl|\sum_{n\le x} \chi(n) 
f(n) \biggr|^{2k}.
$$
\endproclaim

\demo{Proof} 
Write $\Delta=\Delta_f(x,q)$, and define 
$$
\Delta_0:=\sum\limits_{n \le x} \chi_0 (n),\qquad 
 \Delta_1:=\biggl| \sum\limits_{n \le x} \chi_0 (n) f(n)\biggr| \qquad
\text{and} \qquad \Delta_2:=\sum\limits_{n \le x} \chi_0 (n) |f(n)|^2.
$$
Note that the Cauchy-Schwarz inequality gives 
$\Delta_1^2\leq \Delta_0\Delta_2$. 
A straightforward computation, using the orthogonality relations
for characters, gives that 
$$
\Delta_2 =
\frac{1}{\phi (q)} \sum_{\chi (\text{\rm mod}\ q)} 
\left| \sum_{n \le x} \chi (n) f(n)\right|^2 ,
$$ 
and thus 
$$
\left( \Delta_2-\frac{\Delta_1^2}{\phi(q)}\right)^k \le \biggl( \frac{1}{\phi(q)} \sum\Sb\chi \pmod q\\ \chi\ne \chi_0\endSb \biggl| 
\sum_{n \le x} \chi(n) f(n) \biggr|^{2k} \biggr) \biggl(\frac{1}{\phi(q)} 
\sum\Sb\chi \pmod q\\ \chi\ne \chi_0\endSb 1 \biggr)^{k-1} \leq \Delta^{2k} ,
$$
by H{\" o}lder's inequality. This then implies  
$$
\Delta_1^2 \left( 1-\frac{\Delta_0}{\phi(q)}\right) = 
\Delta_1^2 - \frac{\Delta_0\Delta_1^2}{\phi(q)} \leq 
\Delta_0 \left( \Delta_2-\frac{\Delta_1^2}{\phi(q)}\right)
\leq \Delta_0 \Delta^2 .
$$
If $k\geq 2$ then $x\leq \sqrt{q}$ and so $(1-\Delta_0/\phi(q))^k=1+o(1)$,
which, combined with the line above, implies that  $\Delta_1^{2k} \lesssim\Delta_0^k \Delta^{2k}$. Therefore
$$
\frac{1}{\phi (q)} \sum_{\chi (\text{\rm mod}\ q)} 
\left| \sum_{n \le x} \chi (n)f(n) \right|^{2k} 
\leq
\frac{\Delta_1^{2k}}{\phi(q)} + \frac{(\phi
(q)-1)\Delta^{2k}}{\phi(q)} \lesssim \Delta^{2k}\biggl(
\frac{\Delta_0^{k}}{\phi(q)} + 1\biggr). \tag{2.1}
$$
By the small sieve we know that for $x \ge \log q$ 
$$
\Delta_0 = \sum\Sb n \le x \\ (n,q)=1\endSb 1 \le c\frac{\phi (q)}{q} x 
$$
for some absolute constant $c >0$.  Hence 
$$
\Delta_0^k \le x^{k-1} \Delta_0 \le c\frac{\phi(q)}{q} x^k \le c \phi(q)
$$
and the Proposition follows upon inserting this estimate in (2.1).

\enddemo

We cannot expect to get good lower bounds for $\Delta_f(x,q)$ for all 
arithmetic 
functions $f$, since there may be a good deal of cancellation in 
determining the sum $d_{k,f}(n,x)$,
making $\sum_n |d_{k,f}(n,x)|^2$ small.  We shall focus on a large class
${\Cal F}$ of arithmetic functions defined as follows:  $f\in {\Cal F}$ if
$f(n) = g(n) h(n)$ where $g$ is a multiplicative function with $|g(n)|=1$ 
for all $n$, and $h(n)\ge 0$ for all $n$.  Note that ${\Cal F}$ includes
$\mu(n)$ (the M{\" o}bius function), $\omega(n)$ (the number of distinct 
prime divisors of $n$), $d(n)$ (the divisor function), $n^{it}$ 
(for a real number $t$) among others.  

\proclaim{Lemma 2.3}  Suppose $f$ and $g$ are arithmetic functions with 
$f(n) \ge g(n) \ge 0$ for all $n$.  Then for all integers $k\ge 1$
$$
{\Bbb E}\biggl(\biggl|\sum_{n\le x} X_n f(n)\biggr|^{2k} \biggr) 
\ge {\Bbb E}\biggl(\biggl|\sum_{n\le x} X_n g(n) \biggr|^{2k} \biggr).
$$
If $f\in {\Cal F}$ with $|f(n)| \ge \theta$ for all square-free $n$ then
$$
{\Bbb E}\biggl( \biggl|\sum_{n\le x} X_n f(n)\biggr|^{2k}\biggr) 
\ge \theta^{2k} \sum_{N\le x^k} \mu(N)^2 d_{k}(N,x)^2.
$$
\endproclaim
\demo{Proof} If $f(n) \ge g(n) \ge 0$ then $d_{k,f}(n,x) \ge d_{k,g}(n,x)$ 
and so the first assertion follows from Lemma 2.1.  If $f\in {\Cal F}$ and 
$|f(n)|\ge \theta $ for all $n$, then for squarefree $N$ we have
$|d_{k,f}(N,x)| \ge \theta^{k} d_{k}(N,x)$ and so the second statement follows
from the first part of the lemma.
\enddemo

In \S 3 we collect together several results from multiplicative number 
theory; chiefly on smooth numbers (integers not having large prime factors),
and round numbers (integers having many prime factors).  We shall use these 
in \S 4 to estimate the $2k$-th moments of $\sum_{n\le x}^{\flat} X_n$ where 
the flat ``$\flat$'' indicates that the sum is over squarefree $n$ coprime to $q$; and 
in \S 6 to get good estimates for large moments of $\sum_{n\le x} X_n$.  
We show in \S 5 how the estimates of \S 4 lead to the large character sums 
given in Theorems 4 through 7, and Corollaries 1, 2, and 3.  We note that 
these results depend only on the lower bounds for $\sum_{N\le x^k}^\flat 
d_{k}(N,x)^2$ given in Theorems 4.1, and 4.2.  In view of Lemma 2.3 
we may thus generalize these results for $\Delta_f$ when $f\in {\Cal F}$ 
with $|f(n)|\ge 1$.  

\proclaim{Theorems 4-7, Corollaries 1-3 Revisited}  Let $f\in {\Cal F}$ 
be any arithmetic function with $|f(n)| \ge 1$ for all $n$.  Then 
Theorem 4-7 and Corollaries 1-3 all hold for $\Delta_f(x,q)$ in place of 
$\Delta (x,q)$.
\endproclaim

In \S 7 we derive Theorems 1 and 3 as consequences of the analysis of \S 6.  
In \S 8 we obtain the condtional result Theorem 2.  The case of real characters
(Theorems 9-11) are dealt with in \S 9.  Lastly, Theorem 8, which is a 
consquence of the ``Fourier flip'' $x\to \frac{q}{x}$, is proved in \S 10.

\head 3.  Smooth and round numbers \endhead 

\subhead 3a. Integers with a specified number of prime
factors
\endsubhead

\noindent Estimating $\pi(x,y)$, the number of integers up to $x$ with exactly
$y$ distinct prime factors, has long been a central topic of additive number 
theory. Hardy and Ramanujan [8] established the famous upper bound 
$$
\pi (x,y)\ll \frac{x}{\log x} \frac{ (\log\log x + O(1))^{y-1}}{(y-1)!},
$$
uniformly for all $y$.  However good lower bounds, 
even on the order of magnitude for $\pi(x,y)$, when $y\gg \log \log x$ 
were not known until recently.  In 1984, Pomerance [15]
made an important breakthrough in showing that
$$
\pi (x,y) = \frac{x}{\log x} \frac{L^{y+ O(\frac yL)}}{y!}, 
\qquad \text{where    } L = \log\biggl(\frac{\log x}{y\log y}\biggr), 
\tag{3.1}
$$
in the range 
$$
\log\log x \le y \le \frac{\log x}{3\log\log x}. \tag{3.2} 
$$

Pomerance only claimed to have proved this result in the narrower range 
with $y\geq (\log\log x)^2$. However he gives a slightly worse error 
term in one place in his proof than is necessary, 
with the resulting loss in the range of applicability. 
This mistake is corrected in the proof of Theorem 3.1 below; 
taking $m=1$ there implies the lower bound in (3.1).  The upper bound in the
missing range follows from Hardy and Ramanujan's result.

Although it appears that we have imposed some rather severe 
extra restrictions, it turns out that we can obtain the following result 
with minor modifications to Pomerance's proof.  Here $\sum^{\flat}$ indicates 
that the sum is over squarefree arguments.

\proclaim{Theorem 3.1} Given integers $x, y$ and $m$, let $z=\max  
(y^2, \omega (m))$. If (3.2) holds and, in addition,
$$
y^2 \le z \le x^{\frac{2}{3y}} \tag{3.3}
$$
then
$$
\mathop{{\sum}^{\flat}}\Sb n\le x,\ \omega(n)=y \\ (n,m)=1 \endSb 1 
\ge \frac{x}{\log x} \frac{L^{y+ O(\frac yL)}}{y!}, \qquad \text{where    }
 L = L(x,y,z): = \log \left( \frac{\log x}{y \log \sqrt{z}} \right).
\tag{3.4}
$$
\endproclaim

To prove this Theorem we require the following lemma.

\proclaim{Lemma 3.2}  Let $I$ be any interval, and let $s\ge 2$ be an 
integer.  Then
$$
\biggl(\sum_{p\in I} \frac 1p \biggr)^s - \frac{s(s-1)}{2} 
\biggl(\sum_{p \in I} \frac{1}{p^2}\biggr)\biggl(\sum_{p\in I} \frac 1p 
\biggr)^{s-2} \le \sum\Sb p_1,p_2,\ldots,p_s \in I\\ p_i \text{ distinct } 
\endSb \frac{1}{p_1 \ldots p_s } \le \biggl(\sum_{p\in I} \frac 1p\biggr)^s.
$$
\endproclaim
\demo{Proof} The upper bound is immediate, and the lower bound 
follows by induction on $s$, after noting that 
$$
\align
\sum\Sb p_1,\ldots,p_s\in I\\ p_i \text{ distinct }\endSb \frac{1}{p_1\ldots 
p_s} &= \sum\Sb p_1\ldots p_{s-1} \\ p_i \text{ distinct } \endSb 
\frac{1}{p_1 \ldots p_{s-1}} \biggl( \sum_{p\in I} \frac 1p - \frac{1}{p_1} 
-\ldots -\frac{1}{p_{s-1}}\biggr) 
\\
&\ge \biggl(\sum_{p\in I} \frac 1p\biggr) 
\sum\Sb p_1\ldots p_{s-1} \\ p_i \text{ distinct } \endSb 
\frac{1}{p_1 \ldots p_{s-1}} - (s-1) \sum_{p\in I} \frac 1{p^2} 
\sum\Sb p_1\ldots p_{s-2} \\ p_i \text{ distinct } \endSb 
\frac{1}{p_1 \ldots p_{s-2}}.\\
\endalign
$$

\enddemo

\demo{Proof of Theorem 3.1} If $m$ is an integer with exactly $k$ 
distinct prime factors, and $p_k$ is the $k$th smallest prime, then
$$
\mathop{{\sum}^{\flat}}\Sb n\le x ,\ \omega (n) = y \\(n, m) = 1\\ \endSb
1 \ge
\mathop{{\sum}^{\flat}}
\Sb n\le x,\ \omega (n) = y\\ p|n \implies p > p_k\\ 
\endSb 1.
\tag{3.5}
$$
This is evident from noting that if $q_1, \cdots , q_l$ are the 
distinct prime factors of $m$ that are $> p_k$, and $r_1 , \cdots
r_l$ are the primes $\le p_k$ that do not divide $m$, then each
integer $q_1q_2 \cdots q_lt$ counted in the sum on
the right side of (3.5), corresponds to a distinct integer
$r_1r_2\cdots r_lt$ counted in the sum on the left side. 

Note that $L\ge \log 3$, and put $s=[\frac{y-1}{L+20}]$ and 
$J=[\log (L+20)]-2$.  We define the intervals $I_{-1}= (z,x^{\frac{2}{ey}}]$,
and (for $0\le j\le J-1$) $I_j=(x^{\frac{2e^{j-1}}{y}},x^{\frac{2e^j}{y}}]$.
We get a lower bound on the right side of (3.5) by counting only those integers
$n$ of the form $n=n_{-1} n_0\ldots n_{J-1} p$ where $n_{-1}$ consists 
of exactly $y-1-sJ$ distinct primes from $I_{-1}$, and (for $0\le j\le J-1$)
$n_j$ consists of exactly $s$ distinct primes from $I_j$, and $
x^{\frac{2e^{J-1}}{y}}< p \le \frac{x}{n_{-1}n_0\ldots n_{J-1}}$ is prime.
Hence using the prime number theorem 
$$
(3.5) \ge \sum_{n_{-1}, n_0, \ldots, n_{J-1}} \ \ \ \ 
\sum_{x^{\frac{2e^{J-1}}{y}}
\le p \le \frac{x}{n_{-1}n_0\ldots n_{J-1}} } 1 \ge \frac{x}{2\log x} 
\prod_{j=-1}^{J-1} \biggl(\sum_{n_j} \frac 1{n_j}\biggr).
$$

Appealing to Lemma 3.2 we determine that
$$
\sum_{n_{-1}} \frac{1}{n_{-1}} = \frac{(L+O(1))^{y-1-sJ}}{(y-1-sJ)!} \biggl(
1+O\biggl(\frac{y^2}{z\log z}\biggr)\biggr),
$$
and that for $0\le j\le J-1$ 
$$
\sum_{n_j} \frac 1{n_j} = \frac{(1+O(L^{-1}))^s}{s!} \biggl(1+ O\biggl(
\frac{s^2}{x^{\frac{2e^{j-1}}{y}}\log z}\biggr)\biggr).
$$
It follows that 
$$
(3.5) \ge \frac{x}{2\log x} \frac{L^{y-1-sJ}}{(y-1-sJ)! s!^J} e^{O(\frac yL 
+ \frac{1}{\log y})},
$$
and Theorem 3.1 follows upon using Stirling's formula, keeping in mind that 
$y\ge L$ in our range.

\enddemo

\proclaim{Lemma 3.3}  Let $\ell$ be a positive integer, and suppose 
$y\ge 2\ell^2$.  For all $x\ge y^{\ell}$,
$$
\sum\Sb n\le x, \Omega(n)=\ell\\ p|n \implies p>y\endSb  1
\ll \frac{x}{\log y} \frac{(\log \log x +O(1))^{\ell -1}}{(\ell -1)!}.
$$
\endproclaim

\demo{Proof}  Given a squarefree integer $m$ with exactly $j 
(\le \ell)$ distinct prime factors all larger than $y$, 
there are $\le j^{\ell- j}$ integers $n$ with $\Omega(n)=\ell$ and having
exactly the same prime factors as $m$. Moreover if $n\leq x$ then
$m\leq n/y^{\ell-j} \leq x/y^{\ell-j}$. Thus the sum we seek is 
$$
\le \sum_{j=1}^{\ell} j^{\ell- j} \sum\Sb m\le x/y^{\ell-j}\\ \omega(m)=j 
\endSb \mu(m)^2.
$$
By the Hardy-Ramanujan upper bound this is 
$$
\ll \frac{x}{\log y} \sum_{j=1}^{\ell} \frac{j^{\ell-j} }{y^{\ell-j}}
\frac{(\log \log x+O(1))^{j-1}}{(j-1)!} 
\le \frac{x}{\log y} \frac{(\log \log x +O(1))^{\ell-1}}{(\ell-1)!} 
\sum_{j=1}^{\ell}  \frac{(j\ell)^{\ell-j}}{y^{\ell -j}},
$$
and the result follows as $y \ge 2\ell^2$.
\enddemo

\subhead 3b. Smooth numbers \endsubhead

\noindent Given real numbers $x\ge y\ge 1$, and an integer $\ell$ we 
define ${\Cal S}_\ell (x,y)$ to be the set of integers below $x$ having 
exactly $\ell$ prime factors (counted with multiplicity) larger than $y$.  
We denote the cardinality of ${\Cal S}_{\ell}(x,y)$ by $\Psi_{\ell} (x,y)$.  
The case $\ell=0$ gives rise to smooth numbers: that is, integers free of 
large prime factors, and we write ${\Cal S}(x,y)$, $\Psi(x,y)$ in place 
of ${\Cal S}_0(x,y)$, $\Psi_0(x,y)$.  Estimating $\Psi(x,y)$ has been 
the focus of much attention, and we quote below the best results 
known.

\proclaim{Theorem 3.4}  Let $x\ge y\ge 2$ be real numbers and put 
$u=\frac{\log x}{\log y}$.  For any fixed $\epsilon >0$ the 
asymptotic formula
$$
\Psi(x,y) =  x\rho(u) \biggl(1+O\biggl(\frac{\log(u+1)}{\log y}\biggr)\biggr),
\tag{3.6}
$$
holds uniformly in the range $1\le u\le \exp((\log y)^{\frac 35 -\epsilon})$.
The weaker relation
$$
\log \frac{\Psi(x,y)}{x} = \left( 1+ O\left(\exp(-(\log u)^{\frac 35 
-\epsilon})\right)\right)\log \rho(u) 
\tag{3.7}
$$
holds uniformly in the range $1\le u\le y^{1-\epsilon}$.  Lastly, 
as $u\to \infty$
$$
\log \rho(u) = -u \biggl( \log u + \log \log (u+2) - 1+ O\biggl(
\frac{\log \log(u+2)}{\log (u+2)}\biggr)\biggr).
\tag{3.8}
$$
\endproclaim

\demo{Proof}  See Theorems 1.1, 1.2, and Corollary 2.3 of [10].
\enddemo

We next give a bound for $\Psi(x/z,y)$ in terms 
of $\Psi(x,y)$. 

\proclaim{Proposition 3.5}  There is an absolute constant $c$ such that 
for all $1\le z\le x$ and $y\geq 2$, 
$$
\Psi(\tfrac xz, y) \ll (c\log x)^{\frac{\log z}{\log y}}
\frac{\Psi(x,y)}{z}.
$$
\endproclaim

\demo{Proof}  We prove this when $1\le z\le y$; 
the general case follows by repeated application of this result.  
  From Corollary 1.7 of [10] we obtain
$\Psi(\tfrac xz,y) \leq c_1 \Psi(x,y) z^{-\alpha(x/z,y)}$
where $c_1>0$ is some absolute constant, and $\alpha =\alpha(x/z,y)$ is the unique positive solution to 
$\log (x/z) = \sum_{p\le y} \log p/(p^{\alpha}-1)$.  Notice that
$$
\log x \ge \log (x/z) \ge \sum_{n\le y} \frac{\Lambda(n)}{n^{\alpha}} 
\ge \frac{1}{y^{\alpha}} \sum_{n\le y} \Lambda(n) \geq \frac{y}{4y^{\alpha}}.
$$
This shows that $y^{-\alpha}\leq 4(\log x)/y$ so that $z^{-\alpha} \leq 
(4\log x)^{ \frac{\log z}{\log y}}/z$.  The result 
$ \Psi(x/z,y) \leq c_1(4\log x)^{\log z/\log y} \Psi(x,y)/z$ follows for 
 $1\le z\le y$, and repeated applications of this result give
$\Psi(x/z,y) \leq c_1(4c_1\log x)^{\log z/\log y} \Psi(x,y)/z$
in general.
\enddemo

We note here a useful corollary of this result:

\proclaim{Corollary 3.6} Let $0\le \kappa <1$, and let $c$ be as in Proposition
3.5.  Suppose $y \ge e (c\log x)^{\frac{1}{1-\kappa}}$.  Then 
$$
\sum_{n\in {\Cal S}(x,y)} \frac{1}{n^{\kappa}}  \ll  \frac{\log y}
{1-\kappa} \frac{\Psi(x,y)}{x^{\kappa}}, \qquad
\text{and} \qquad \sum_{n\in {\Cal S}(x,y)} \frac{1}{n^{\kappa}} 
\log \frac{x}{n} \ll \frac{\log y}{(1-\kappa)^2} \frac{\Psi(x,y)}{x^{\kappa}}.
$$
If $y \ge (c\log x)^2$ then 
$$
\sum_{n\in {\Cal S}(x,y)} \log \frac xn \ll \Psi(x,y).
$$
\endproclaim

\demo{Proof} By partial summation 
$$
\sum_{n\in {\Cal S}(x,y)} \frac{1}{n^{\kappa}} = \int_{1^{-}}^{x} \frac{1}
{t^{\kappa}} d\Psi(t,y) = \frac{\Psi(x,y)}{x^{\kappa}} + \kappa \int_1^x 
\frac{\Psi(t,y)}{t^{1+\kappa}} dt.
$$
Using Proposition 3.5 the second term above is 
$$
\ll \kappa \frac{\Psi(x,y)}{x} \int_1^x \frac{1}{t^{\kappa}} 
(c\log x)^{\frac{\log (x/t)}{\log y}} dt =  
\kappa \frac{ \Psi(x,y) }{x} \int_1^x 
\frac{x^{\frac{\log (c\log x)}{\log y}}}{t^{\kappa +\frac{\log (c\log x)}
{\log y}}}dt,
$$
and using our hypothesis on $y$ this is 
$$
\le \frac{\kappa \log y}{1-\kappa} \frac{\Psi(x,y)}{x^{\kappa}}.
$$
The first part of the corollary follows.  The other two assertions are
proved similarly.
\enddemo

\proclaim{Lemma 3.7}  Let $x\ge y \ge (\log x)^{1+\epsilon}$, and put 
$u = \frac {\log x}{\log y}$.  Then
$$
\frac{\Psi(x,y\log y)}{\Psi(x,y)} 
= \exp\biggl( u\frac{\log \log y}{\log (y\log y)} (\log u+O(\log \log (u+2)))
\biggr).
$$
\endproclaim
\demo{Proof}  By Lemma 2.2 and Corollary 2.4 of [10] we get
$$
\frac{\rho(\frac{\log x}{\log (y\log y)})}{\rho(\frac{\log x}{\log y})} 
= \exp\biggl( u\frac{\log \log y}{\log (y\log y)} (\log u+O(\log \log (u+2)))
\biggr).
$$
The lemma follows upon combining this with (3.6) when 
$u\le \exp((\log y)^{\frac 35-\epsilon})$, and (3.7)
for larger $u$.
\enddemo

\proclaim{Lemma 3.8}  Suppose $y\ge (\log x)^{\frac 32}$ and that 
$x\ge z\ge x y^{-\frac 13}$.  Then 
$$
\Psi(x+z,y)-\Psi(x,y) \gg z \frac{\Psi(x,y)}{x}.
$$
\endproclaim 

This follows from Theorems 5.1 and 5.2 of [10].
The next result is an immediate consequence of Corollary 2.4 of [10].

\proclaim{Lemma 3.9} $\rho(u-v)\asymp \rho(u)$ if $|v| \ll 1/\log 2u$,
for $u, u-v\geq 1$. \endproclaim 

\head 4.  The $2k$-th moment of $\sum_{n\le x}^{\flat} X_n$
\endhead

\noindent In this section we prove upper and lower bounds on 
the $2k$-th moment of $\sum_{n\le x}^{\flat} X_n$, where (throughout 
this section) the $\flat$ 
indicates that the sum is over squarefree $n$ coprime to $q$.  These 
bounds will be useful in deducing many of our 
large character sums results.  

\proclaim{Theorem 4.1}  Let $k\ge 1$ be an integer, and 
put $K=\max(k,\omega(q))$.  
Uniformly for all $x \ge K^{ek}$ we have
$$
{\Bbb E}\biggl( \biggl|\mathop{{\sum_{n\le x}}^{\flat}} X_n\biggr|^{2k} 
\biggr)^{\frac1{2k}} \le 
{\Bbb E}\biggl( \biggl|\sum_{n\le x} X_n\biggr|^{2k}\biggr)^{\frac{1}{2k}} 
\le x^{\frac 12} \biggl(\frac{\log x}{k}\biggr)^{\frac{(k-1)^2}{2k}} e^{O(k)},
\tag{4.1}
$$
and 
$$
\align
{\Bbb E}\biggl( \biggl|\mathop{{\sum_{n\le x}}^{\flat}
} X_n\biggr|^{2k} \biggr)^{\frac
1{2k}} &\ge \biggl(\mathop{{\sum_{N\le x^k}}^{\flat}} d_k(N,x)^2 
\biggr)^{\frac{1}{2k}}
\\
&\ge 
\frac{ x^{\frac 12}}{(\log x)^{1-\frac{1}{2k}}} \biggl(\frac{\log x}{k\log K}
\biggr)^{\frac{k}{2}} \biggl(\log \biggl(\frac{\log x}{k\log K}\biggr)
\biggr)^{O(k)}. \tag{4.2}
\\
\endalign
$$
\endproclaim

\demo{Proof of the upper bound (4.1)}  Observe that
$$
\align
{\Bbb E}\biggl(\biggl|\mathop{{\sum_{n\le x}}} X_n \biggr|^{2k}
\biggr) &= \sum_{N\le x^k} d_k(N,x)^2 
\\
&= 
\#\{ b_1,b_2, \ldots, b_k , B_1, B_2, \ldots, B_k \le x: b_1 \ldots b_k = B_1
\ldots B_k \}.\\
\endalign
$$
To each solution above we associate a $k\times k$ ``g.c.d.-matrix'' of 
integers $A=(a_{i,j})$ defined as follows:  Put $a_{1,1} = (b_1,B_1)$, 
and then define (using induction on $i+j$) 
$$
a_{i,j} = \biggl( \frac{b_i}{\prod_{\ell < j} a_{i,\ell}}, \frac{B_j}
{\prod_{\ell < i} a_{\ell,j}}\biggr),
$$
so that, for $1\le i\le k$,
$$
b_i = \prod_{\ell =1}^{k} a_{i,\ell}, \qquad \text{and } \qquad 
B_i = \prod_{\ell =1}^{k} a_{\ell,i}.
$$

We will bound the number of $k\times k$ integer matrices 
$A=(a_{i,j})$ with all row and column products $\prod_{\ell =1}^{k} 
a_{i,\ell}$,  $\prod_{\ell= 1}^{k} a_{\ell, i}$ less than $x$,
which thus implies an upper bound in our original problem.  
The number of choices for $a_{k,k}$ is 
$$
\le \min\biggl(\frac{x}{\prod_{i=1}^{k-1} a_{i,k}},\frac{x}{\prod_{i=1}^{k-1} 
a_{k,i}}\biggr) \le \frac{x}{\prod_{i=1}^{k-1} (a_{i,k}a_{k,i})^{\frac 12}}.
$$
Next we sum over the possibilities for $a_{i,k}$, $a_{k,i}$ ($1\le i\le k-1$).
Notice that $a_{i,k} \le x/\prod_{j=1}^{k-1} a_{i,j}$ and $a_{k,i} 
\le x/\prod_{j=1}^{k-1} a_{j,i}$, and so 
$$
\sum_{a_{i,k}} \frac{1}{\sqrt{a_{i,k}}} 
\le \frac {2\sqrt{x}}{ \prod_{j=1}^{k-1}
a_{i,j}^{\frac 12}}, \qquad \qquad
\sum_{a_{k,i}} \frac{1}{\sqrt{a_{k,i}}}
 \le \frac {2\sqrt{x}}{ \prod_{j=1}^{k-1}
a_{j,i}^{\frac12}}.
$$
Thus given $a_{i,j}$ ($1\le i,j\le k-1$), the number of possibilities 
for the last row and column of $A$ is 
$$
\le \frac{2^{2k-2} x^k}{\prod_{1\le i,j\le k-1} a_{i,j}}.
$$

We now sum this over all the possibilities for $a_{i,j}$ ($1\le i,j \le k-1$).
Keeping in mind that $\prod_{j=1}^{k-1} a_{i,j} \le x$ for all $1\le i\le k-1$,
we see that this is
$$
\le 2^{2k-2} x^k \prod_{i=1}^{k-1} \biggl(\sum_{a_{i,1} ... a_{i,k-1}\le x} 
\frac{1}{a_{i,1}\ldots a_{i,k-1}} \biggr) = 2^{2k-2} x^k \biggl(\sum_{n\le x} 
\frac{d_{k-1}(n)}{n}\biggr)^{k-1}.
$$

Now for any $\alpha >0$ 
$$
\sum_{n\le x} \frac{d_{k-1}(n)}{n} \le x^{\alpha} \sum_{n\le x} 
\frac{d_{k-1}(n)}{n^{1+\alpha}} \le x^{\alpha} \zeta(1+\alpha)^{k-1} = 
x^{\alpha} \biggl(\frac 1\alpha +O(1) \biggr)^{k-1}.
$$
Choosing (optimally) $\alpha = k/\log x$ we obtain (since $k\le \log x$)
$$
\sum_{n\le x} \frac{d_{k-1}(n)}{n} \le \biggl(\frac{\log x}{k}\biggr)^{k-1} 
e^{O(k)}.
$$

To sum up, we have shown that 
$$
{\Bbb E}\biggl(\biggl|\mathop{{\sum_{n\le x}}} X_n \biggr|^{2k}
\biggr) \le 2^{2k-2} x^k \biggl(\sum_{n\le x} \frac{d_{k-1}(n)}{n}\biggr)^{k-1}
\le x^k \biggl(\frac{\log x}{k}\biggr)^{(k-1)^2} e^{O(k^2)},
$$
and (4.1) follows.

\enddemo

\demo{Proof of the lower bound (4.2)}   We bound 
$$
{\Bbb E}\biggl(\biggl|\mathop{{\sum_{n\le x}}^{\flat}} 
X_n \biggr|^{2k}\biggr) 
\ge \mathop{{\sum_{N\le x^k}}^{\flat}} d_{k}(N,x)^2 
$$
by focussing only on special values of $N$ for which we expect $d_k(N,x)$ 
to be large.  Specifically, we let $y$ denote an integer parameter to be 
chosen later, and consider only those $N\le x^k$ which are square-free,
coprime to $q$, 
and have $ky$ distinct prime factors.  Using Cauchy's inequality, we
find that 
$$
\mathop{{\sum_{N\le x^k}}^{\flat}} d_{k}(N,x)^2 
\ge 
\biggl(\mathop{{\sum\Sb N\le x^k \\ \omega(N)=ky\endSb}^{\!\!\!\!\!\!\!\flat}}
 d_k(N,x)
\biggr)^2 \frac{1}{\pi(x^k,ky)}. \tag{4.3}
$$

We shall choose $y=[kL_0]$ where $L_0 =\log (\frac{\log x}{k\log K})$. 
Using (3.1) (after checking that the constraint (3.2) is met) we find that 
$$
\pi(x^k,ky) = \frac{x^k}{k\log x} \frac{L^{ky+O(\frac{ky}{L})}}{(ky)!} 
\qquad \text{ where    } L = \log \biggl(\frac{\log x}{y\log (ky)}\biggr).
$$
Since $L_0 -L = \log \frac yk + \log \frac{\log (ky)}{\log K} \ll \log L_0$ 
we conclude that 
$$
\pi(x^k,ky) = \frac{x^k}{k\log x} \frac{L_0^{ky}}{(ky)!} 
\exp\biggl(O\biggl(ky \frac{
\log L_0}{L_0} \biggr)\biggr). \tag{4.4}
$$

Next observe that 
$$
\mathop{{\sum\Sb N\le x^k\\ \omega(N)=ky\endSb}^{\!\!\!\!\!\!\flat}}
\ d_k(N,x) \ge \mathop{{\sum\Sb m_1,\ldots,m_k\le x\\ \omega(m_i)=
y\endSb}^{\!\!\!\!\!\!\!\!\!\!\!\!*} } \ \ \ 1
$$
where the $*$ indicates that the sum is over squarefree $m_1$ coprime 
to $q$, and pairwise coprime.  We deduce from (3.5) that this is 
$$
\ge  \biggl( \sum\Sb n\le x,  \ \ \omega(n)=y \\ p|n \implies
p>p_{\ell}\endSb \mu(n)^2 \biggr)^k \qquad \text{where     } \ell=(k-1)y 
+ \omega(q).
$$
Now we use Theorem 3.1 to bound this quantity.  Our assumption that $k\log K
\le e^{-1} \log x$ ensures that the criteria (3.2) and (3.3) are met.  Hence,
with $z=\max(\ell,y^2)$, 
$$
 \sum\Sb n\le x,  \ \ \omega(n)=y \\ p|n \implies
p>p_{\ell}\endSb \mu(n)^2 \ge \frac{x}{\log x} \frac{L_1^{y+O(\frac y{L_1})}}
{y!} \qquad \text{where     } L_1 =\log \biggl(\frac{\log x}{y\log \sqrt{z}}
\biggr).
$$
Since $L_0-L_1\ll \log L_0$ we obtain 
$$
\mathop{{\sum\Sb N\le x^k\\ \omega(N)=ky\endSb}^{\!\!\!\!\!\!\!\flat}}
\ \  d_k(N,x) \ge \frac{x^k}{(\log x)^k} \frac{L_0^{ky}}{y!^k} 
\exp\biggl(O\biggl(ky \frac{\log L_0}{L_0}\biggr)\biggr).
\tag{4.5}
$$

Using (4.4) and (4.5) in (4.3) we deduce 
$$
\biggl(\mathop{{\sum_{N\le x^k}}^{\flat}} 
\ d_k(N,x)^2 \biggr)^{\frac 1{2k}} 
\ge \frac{x^{\frac 12}}{(\log x)^{1-\frac{1}{2k}} }\frac{L_0^{\frac y2} 
(ky)!^{\frac{1}{2k}}}{y!}\exp\biggl(O\biggl(\frac{y\log L_0}{L_0}
\biggr)\biggr),
$$
and the lower bound (4.2) follows upon using Stirling's formula, and 
recalling the definitions of $y$ and $L_0$.

\enddemo

Next we give lower bounds on the $2k$-th moment of $\sum_{n\le x}^{\flat}
X_n$ when $x$ is small (roughly, $x =K^A$ for some integer A).

\proclaim{Theorem 4.2}  Let $k$, $A$ be positive integers, and put 
$K=\max(k,\omega(q))$.  For all $x \ge (4(Ak+K)\log (Ak+K))^A$ 
we have, uniformly,
$$
{\Bbb E}\biggl(\biggl|\mathop{{\sum_{n\le x}}^{\flat}} X_n\biggr|^{2k}
\biggr)^{\frac 1{2k}}
\ge \biggl( \mathop{{\sum_{N\le x^k}}^{\!\!\flat}} \ d_{k}(N,x)^2
\biggr)^{\frac{1}{2k}}
\gg    x^{\frac 12} \biggl(\frac{k}{\log x}\biggr)^{\frac A2}.
\tag{4.6} 
$$
\endproclaim

\demo{Proof} Plainly
$$
{\Bbb E}\biggl(\biggl|\mathop{{\sum_{n\le x}}^{\flat}} 
X_n\biggr|^{2k}\biggr)
= \mathop{{\sum\Sb N\le x^k\endSb}^{\flat}} d_k(N,x)^2 
\ge \mathop{{\sum\Sb N\le x^{k} \\ (N,q)=1\endSb}^{\!\!\!\!\!\!*}} d_k(N,x)^2
$$
where the $*$ indicates that we sum over only those $N$ that are 
square-free and composed of
exactly $Ak$ prime factors, all less than $x^{1/A}$.  Note that 
for such $N$, $d_k(N,x)$ is at least the number of $k$-tuples $m_1$, $\ldots$, 
$m_k$ whose product is $N$, where each $m_i$ is the product of 
exactly $A$ primes. Thus $d_k(N,x)\ge (Ak)!/A!^k$, and so 
$$
\align
\mathop{{\sum_{N\le x^k}}^{\!\!\flat}} \ d_{k}(N,x)^2
&\ge \frac{(Ak)!}{(A!)^k} \mathop{{\sum\Sb N\le x^k \\ (N,q)=1\endSb}^{\!\!\!\!
\!\!*}} 
d_k(N,x)
\ge \frac{(Ak)!}{(A!)^k}\frac{1}{(A!)^k}
\sum\Sb p_1, \ldots, p_{Ak}\le x^{\frac 1A} 
\\ p_i\neq p_j, p_i \nmid q\endSb 1
\\
&\ge\frac{(Ak)!}{(A!)^{2k}}
 \prod_{j=1}^{Ak} \biggl(\pi(x^{\frac 1A})-\sum\Sb p|q\\ p\le 
x^{\frac 1A} \endSb 1  - j+1\biggr).\tag{4.7}\\
\endalign
$$ 
By the prime number theorem, and our lower bound for $x$, we get 
$$
\pi(x^{\frac 1A})-\sum\Sb p|q \\ p\le x^{\frac 1A}\endSb 1 - Ak \ge 
\pi(x^{\frac 1A})- Ak -K  
\ge  \frac{Ax^{\frac 1A}}{\log x} - Ak -K \ge \frac{Ax^
{\frac 1A}}{2\log x}.
$$
Using this, and Stirling's formula, in (4.7) we get Theorem 4.2.

\enddemo

\head 5.  Applications to large character sums \endhead

\noindent In this section, we use Theorems 4.1 and 4.2 to deduce many
of our results on large character sums.  We split these results in two 
parts:  when $\log \log x \le (\frac 12 +o(1)) \log \log q$ where 
we use Theorem 4.2, and when $\log \log x \ge (\frac 12+o(1))\log \log q$ 
where Theorem 4.1 is most useful.

\subhead 5a.  Large character sums when
$\log \log x\leq (\frac 12+o(1)) \log\log q$
\endsubhead

\demo{Proof of Theorem 4}  Recall that $x=(10\log q)^B$ for some $B\ge 1$.  
We take $k=[\frac{\log q}{\log x}]$ and $A=[B]$.  Notice that $K
=\max(k,\omega(q)) \le (1+o(1)) \frac{\log q}{\log \log q}$, and so $kA+K \le 
\frac{5}{2}\frac{\log q}{\log \log q}$.  We check now that the condition of 
Theorem 4.2 is met, and so 
$$
\Delta \gg x^{\frac12}\biggl(\frac{k}{\log x}\biggr)^{\frac{[B]}{2}}
\ge \frac{x^{\frac12 +\frac{[B]}{2B}}}{(4\log x)^{[B]}}.
$$
This gives the portion of Theorem 4 not having any 
restriction on $q$.

For our next application we suppose that $\omega(q)\le (\log q)^{\frac{B}{B+1}
-\epsilon}$.  Here, we take $A=[B]+1$, and $k=[x^{\frac 1A}/(10\log x)]$.  Our
bound on $\omega(q)$ ensures that the condition of Theorem 4.2 is met, and so 
$$
\Delta \gg x^{1/2}\biggl(\frac{k}{\log x}\biggr)^{\frac{[B]+1}{2}}
\ge \frac{x}{(4\log x)^{[B]+1}}.
$$
This gives the second part of Theorem 4.  Corollaries 1 and 2 are 
immediate consequences.  
\enddemo

\subhead 5b.  Large character sums when $\log \log x \ge (\frac 12 +o(1)) 
\log \log q$ \endsubhead

\noindent  \demo{Proof of Theorem 5}  We take $k=[\frac{c}{\eta} \sqrt{\frac
{\log q}{\log \log q}}]$ for a fixed, but sufficiently
 small positive constant $c$.  Since $K\le \log q$, one can verify that 
the condition $x\ge K^{ek}$ of Theorem 4.1 is met.  Hence by (4.2) we
get
$$
\Delta \gg \frac{x^{\frac 12}}{(\log x)^{1-\frac 1{2k}} } 
\biggl(\frac{\log x}{k\log\log q}\biggr)^{\frac k2} \biggl(\log \frac{\log x}
{k\log \log q}\biggr)^{O(k)} 
= \frac{x^{\frac 12}}{(\log x)^{1-\frac 1{2k}}} \biggl(\frac{\eta \tau}{c}
\biggr)^{\frac{k}{2}} \biggl(\log \frac{\eta \tau}{c}\biggr)^{O(k)}.\!
$$
The result follows if $c$ is sufficiently small.
\enddemo

\demo{Proofs of Theorems 6 and 7} Both these results follow upon using 
(4.2) with $k=[\frac{\log q}{\log x}]$: the hypotheses in the Theorems 
ensure that $x\ge (\log q)^{ek} \ge K^{ek}$.  

\enddemo
\head 6.  The $2k$-th moment of $\sum_{n\le x} X_n$
\endhead

\noindent Here we explore more finely the $2k$-th moment of $\sum_{n\le x} 
X_n$.  Put
$$
\Psi_\ell (x,y;X_n) = \sum\Sb n\in {\Cal S}_{\ell}(x,y)\endSb X_n.
$$
Our aim in this section is to show that $\sum_{n\le x} X_n$ behaves 
like $\Psi_0(x,y;X_n)$ most of the time, for an appropriately chosen $y$.

\proclaim{Theorem 6.1}  Suppose $k\ge 2$ is an integer, and that 
$y\ge C\log^2 x$ for a large absolute constant $C$.  Then
$$
{\Bbb E}\biggl( \biggl|\sum_{n\le x} X_n -\Psi_0(x,y;X_n)\biggr|^{2k}\biggr)
^{\frac 1{2k}} \ll \Psi(x,y) \biggl(\frac{k\log y\log^2 x}{y}\biggr)^{\frac 12}
\exp\biggl(O\biggl(\frac{k \log^2 x \log\log x }{y}\biggr)\biggr).
$$
\endproclaim

\demo{Proof} Put 
$u=\frac{\log x}{\log y}$.  By Minkowski's inequality 
$$
\align
{\Bbb E}\biggl( \biggl|\sum_{n\le x} X_n -\Psi_0(x,y;X_n)\biggr|^{2k}\biggr)
^{\frac 1{2k}} 
&= {\Bbb E}\biggl( \left|\sum_{\ell =1}^{[u]}\Psi_{\ell}(x,y;X_n)\right|^{2k}\biggr)^{\frac{1}{2k}} \\
&\le \sum_{\ell =1}^{[u]} {\Bbb E}\left(|\Psi_{\ell}(x,y;X_n)|^{2k}\right)
^{\frac{1}{2k}}.
\tag{6.1}
\\
\endalign
$$
Observe that
$$
{\Bbb E}(|\Psi_{\ell}(x,y;X_n)|^{2k}) = \sum\Sb m_1 ... m_k = m_1^{\prime} ...
m_k^{\prime} \\ m_i, m_i^{\prime} \in {\Cal S}(x/y^\ell, y) \endSb 
\ \ \ 
\sum\Sb n_1... n_k = n_1^{\prime} ... n_k^{\prime} \\ n_i\le x/m_i, 
n_i^{\prime} \le x/m_i^{\prime} \\ \Omega(n_i)=\Omega(n_i^{\prime}) =\ell \\
p|n_i, n_i^{\prime} \implies p>y\endSb 1.
$$
Now, given $N=n_1... n_k$, the number of factorizations $N = n_1^{\prime} ... n_k^{\prime}$ with each $\Omega(n_i^{\prime}) =\ell$ is $\leq (k\ell)!/\ell!^k$, and so the inner sum over $n_i$, $n_i^{\prime}$ is 
$$
\leq \sum\Sb n_1, \ldots, n_k\\ n_i\le x/m_i \\ \Omega(n_i)=\ell \\ p|n_i \implies
p>y\endSb \frac{(k\ell)!}{\ell!^k} \ll \frac{k^{k\ell }}{\ell^{(k-1)/2}} \prod_{i=1}^{k} 
\sum\Sb n_i\le x/m_i\\ \Omega(n_i)=\ell \\ p|n_i\implies p>y\endSb 1.
$$
Using Lemma 3.3 (note that $y \ge C \log^2 x\ge 2\ell^2$)
this is
$$
\align
&\le \frac{x^k k^{k\ell}}{m_1 \ldots m_k} \frac{(\log \log x+O(1))^{k(\ell-1)}}
{(\ell-1)!^k} \left( \frac{c}{\log y} \right)^k \frac{1}{\ell^{(k-1)/2}} \\
&\ll \frac{x^k}{m_1 \ldots m_k} \frac{ (2k\log\log x)^{k\ell} } {(\ell! \log y\log \log x)^k} .\\
\endalign
$$

Now
$$
\sum\Sb m_1\ldots m_k = m_1'\ldots m_k'\\ 
m_i,m_i' \in {\Cal S} (x/y^\ell,y) \endSb \frac{1}{m_1\ldots m_k} =
{\Bbb E}\biggl(\biggl|\sum_{n\in {\Cal S}
(x/y^{\ell},y)} \frac{X_n}{\sqrt{n}} \biggr|^{2k}\biggr) \leq 
\biggl(\sum_{n\in {\Cal S}
(x/y^{\ell},y)} \frac{1}{\sqrt{n}} \biggr)^{2k},
$$
and, by Corollary 3.6 and  Proposition 3.5
$$
\sum_{n \in {\Cal S}(x/y^{\ell},y)}\frac{1}{\sqrt{n}} 
\ll \log y \frac{\Psi(x/y^{\ell},y)}{(x/y^{\ell})^{\frac 12}} 
\ll \log y (c\log x)^{\ell}  \frac{\Psi(x,y)}{\sqrt{xy^{\ell}}}.
$$
Therefore, combining the bounds above, we get
$$
{\Bbb E}\left(|\Psi_{\ell}(x,y;X_n)|^{2k}\right)^{\frac{1}{2k}} \ll
\Psi(x,y)  \left( \frac{\log y}{\ell! \log\log x} \right)^{1/2}
 \left(  \frac{c k \log^2 x \log\log x}{y} \right)^{\ell/2}
$$
for some constant $c>0$. Therefore, substituting this into (6.1), we get 
$$
\align
{\Bbb E}\biggl( \biggl|\sum_{n\le x} X_n -\Psi_0(x,y;X_n)\biggr|^{2k}\biggr)
^{\frac 1{2k}}
&\ll  \Psi(x,y) \left( \frac{\log y}{\log\log x}\right)^{1/2} \sum_{\ell=1}^{[u]} \frac{1}{\ell!^{1/2}} \left(  
\frac{c k \log^2 x \log\log x}{y} \right)^{\ell/2}\\
&\ll \Psi(x,y) \left( \frac{ k \log^2 x \log y}{y} \right)^\frac12
\exp\biggl( O\biggl(\frac{k \log^2 x \log \log x}{y}\biggr)\biggr),\\
\endalign
$$
since $\sum_{j=0}^{\infty} \xi^{\frac j2}/j!^{\frac 12} \ll e^{\xi}$ for all
$\xi \ge 0$.  This proves the theorem.

\enddemo

We now derive a good lower bound for the $2k$-th moment of $\sum_{n\le x} X_n$.
This is a considerable refinement of Theorem 4.2, in the case that $q=1$.  

\proclaim{Theorem 6.2}  Let $k\ge 2$ be an integer.  Then for all 
$y\geq 2$ we have 
$$
{\Bbb E}\biggl(\biggl|\sum_{n\le x}X_n
\biggr|^{2k}\biggr)^{\frac 1{2k}} \ge \Psi(x,y) 
\exp\biggl(-\frac{2y\log \log x}{k\log y}
+O\biggl(\frac{1}{\log x}\biggr)\biggr).
$$
\endproclaim
\demo{Proof}  Using Lemma 2.3 with $f(n)=1$, and $g(n)=$ the characteristic
function of ${\Cal S}(x,y)$ we have
$$
{\Bbb E}\biggl(\biggl|\sum_{n\le x}X_n
\biggr|^{2k}\biggr) \ge {\Bbb E}\biggl(\biggl|\sum_{n\in {\Cal S}(x,y)}X_n
\biggr|^{2k}\biggr). 
$$
We bound the right side above by picking only those $X_n$ for 
which $|\arg (X_p)| \le \pi (\log x)^{-2}$ for all $p\le y$ (where $\arg$ is defined to 
lie between $-\pi$ and $\pi$).  
The probability of this happening is clearly $(\log x)^{-2\pi(y)}$.  For 
such a choice of $X_p$'s note that 
$$
\biggl|\sum_{n\in {\Cal S}(x,y)} X_n - \Psi(x,y)\biggr| \le  
\sum_{n\in {\Cal S}(x,y)} |X_n-1| \ll 
\sum_{n\in {\Cal S}(x,y)} \frac{\Omega(n)}{\log^2 x} \ll \frac{\Psi(x,y)}
{\log x}.
$$
Hence 
$$
{\Bbb E}\biggl(\biggl|\sum_{n\in {\Cal S}(x,y)}X_n
\biggr|^{2k}\biggr) \ge \Psi(x,y)^{2k} \exp\biggl(-2\pi(y)\log \log x +
O\biggl(\frac{k}{\log x}\biggr)\biggr),
$$
and the result follows.

\enddemo

Combining Theorems 6.1 and 6.2 we get good upper and lower estimates
for large moments of $\sum_{n\le x} X_n$; and in fact, we get an asymptotic
formula for very large $k$.  

\proclaim{Corollary 6.3} If $k\ge C\log x$ is an integer then
$$
\Psi(x,k\log^5 k \log x) \biggl(1+O\biggl(\frac{1}{\log x}\biggr)\biggr)
\ge {\Bbb E}\biggl(\biggl|\sum_{n\le x}X_n
\biggr|^{2k}\biggr)^{\frac 1{2k}} \ge \Psi(x,\tfrac{ k\log x}{\log^5 k}) 
\biggl(1+O\biggl(\frac{1}{\log x}\biggr)\biggr).
\tag{6.3}
$$
If $\log k/\sqrt{\log x} \log\log x \to \infty$ then 
$$
{\Bbb E}\biggl(\biggl|\sum_{n\le x}X_n
\biggr|^{2k}\biggr)^{\frac 1{2k}} = (1+o(1)) \Psi(x,k).
\tag{6.4}
$$
\endproclaim
\demo{Proof}  From Theorem 6.2 we get 
$$
{\Bbb E}\biggl(\biggl|\sum_{n\le x}X_n
\biggr|^{2k}\biggr)^{\frac 1{2k}} \ge \Psi(x, \tfrac{k \log x}{\log^{3} k})
\exp\biggl(O\biggl(\frac{1}{\log x} + \frac{\log x}{\log^3 k}\biggr)\biggr).
$$
The lower bound of (6.3) now follows upon appealing to Lemma 3.7.  
Using Minkowski's inequality and Theorem 6.1 (with $y= k\log^4 k\log x$)
we get, since $|\Psi_0(x,y;X_n)|\leq \Psi(x,y)$,
$$
\align
{\Bbb E}\biggl(\biggl|\sum_{n\le x}X_n
\biggr|^{2k}\biggr)^{\frac 1{2k}} &\le {\Bbb E}(|\Psi_0(x,y;X_n)|^{2k})
^{\frac{1}{2k}} + O\biggl(\Psi(x,y) \biggl(
\frac{\log x}{\log^3 k}\biggr)^{\frac 12}
\exp\biggl(O\left(\frac{\log x}{\log^3 k}\right)\biggr)\biggr)\\
&\le \Psi(x,k\log^4 k \log x)\exp\biggl(O\left(\frac{\log x}{\log^3 k}
\right)\biggr).
\endalign
$$
The upper bound of (6.3) follows from this and Lemma 3.7.   By Lemma 3.7
we deduce that if $k>\exp ( \sqrt{\log x})$ then
$$
\Psi(x,k(\log k)^{O(1)}) = \Psi(x,k) \exp\biggl( O\biggl(
\log x \frac{(\log\log x)^2}{\log^2  k}\biggr)\biggr).
$$
Therefore (6.4) follows from (6.3).
\enddemo

\head 7.  Implications for character sums: Proofs of Theorems 1 and 3
 \endhead

\noindent Observe that for any integer $k\le \frac {\log q}{\log x}$, 
and any $y$, we have 
$$
\frac{1}{\phi(q) }\sum_{\chi\pmod q} \biggl|\sum_{n\le x}
\chi(n) -\Psi(x,y;\chi)\biggr|^{2k} \le {\Bbb E}\biggl(\biggl|\sum_{n\le x}
X_n - \Psi(x,y;X_n)\biggr|^{2k}\biggr).
$$
Using Theorem 6.1 we deduce that if $y\ge C\log^2 x$ then 
$$
\align
\frac{1}{\phi(q) }\sum_{\chi\pmod q} \biggl|\sum_{n\le x}&
\chi(n) -\Psi(x,y;\chi)\biggr|^{2k} \\
\leq & c^k \Psi(x,y)^{2k} \biggl(
\frac{k\log y \log^2 x}{y}\biggr)^{k} \exp\biggl(O\biggl(
\frac{k^2\log\log x \log^2 x}{y}\biggr)\biggr), \tag{7.1}\\
\endalign
$$
for some constant $c>0$.

\demo{Proof of Theorem 1}  We choose $k=[\frac{\log q}{\log x}]$.  It 
follows from (7.1) that for any $A>1$ there are fewer than $qA^{-2k}$ 
characters $\chi\pmod q$ not satisfying 
$$
\biggl|\sum_{n\le x} \chi(n) -\Psi(x,y;\chi) \biggr| 
\ll A \Psi(x,y) \biggl(\frac{\log q \log x \log y}{y}\biggr)^{\frac 12}
\exp\biggl(O\biggl(\frac{\log q\log x\log \log x}{y}\biggr)\biggr).
$$

Taking $y\ge \log q \log x (\log \log q)^5$, and $A=10$ above, we obtain 
the first assertion of Theorem 1.  

Next, take $y=(\log q +\log^2 x) (\log \log q)^4$, and $A=\exp(\frac{\log x}
{(\log \log q)^2})$.  We deduce that with at most 
$q^{1-\frac{1}{(\log \log q)^2}}$ exceptions
$$
\biggl|\sum_{n\le x} \chi(n)\biggr| \ll \Psi(x,y) \exp\biggl(O\biggl(
\frac{\log x}{(\log \log q)^2}\biggr) \biggr) \ll \Psi (x,y\log y),
$$
using Lemma 3.7.  This gives the second part of Theorem 1.

\enddemo

We now move towards the proof of Theorem 3.  We begin with a lemma 
which may be of independent interest.  

\proclaim{Lemma 7.1}  Let $f(n)$ be any completely multiplicative 
function with $|f(n)|=1$ for all $n$.  Let $2\le x
\le \exp((\log q)^{\frac 12})$,  and let $y=\log q/(\log x (\log \log q)^8)$. 
There are at least $q^{1-\frac{1}{(\log \log q)^2}}$ characters $\chi \pmod q$ 
with 
$$
\sum_{n\in {\Cal S}(x,y)} \chi(n) = \sum\Sb 
n\in{\Cal S}(x,y)\\ (n,q)=1\endSb f(n) + O \biggl(\frac{\Psi(x,y;\chi_0)}
{(\log \log q)^2}\biggr).
$$
\endproclaim

\demo{Proof}  Note that for any integer $k\le \frac{\log q}{\log x}$
$$
\frac{1}{\phi(q)} \sum_{\chi\pmod q} \biggl|\sum\Sb n\in{\Cal S}(x,y)\\
(n,q)=1\endSb \frac{\chi(n) \overline{f(n)}+1}{2} \biggr|^{2k} 
= {\Bbb E} \biggl( \biggl|
\sum\Sb n\in {\Cal S}(x,y)\\ (n,q)=1\endSb \frac{X_n +1}{2}
\biggr|^{2k}\biggr). \tag{7.2}
$$
We give a lower bound for the right side of (7.2) by the argument of Theorem 
6.2.  We pick only those $X_n$ with $|\text{arg}(X_p)| \le \frac{\pi}{\log q}$
for all $p\le y$.  This happens with probability $\ge (\log q)^{-\pi(y)}
\ge \exp(-3y)$, and for such a choice 
$$
\sum\Sb n\in {\Cal S}(x,y)\\ (n,q)=1\endSb \frac{X_n+1}{2} = \Psi(x,y;\chi_0)
+ O\biggl(\sum\Sb n\in{\Cal S}(x,y)\\(n,q)=1\endSb \frac{\Omega(n)}{\log q}
\biggr) = \Psi(x,y;\chi_0)\biggl(1+O\biggl(\frac{\log x}{\log q}\biggr)\biggr).
$$
It follows that 
$$
\frac{1}{\phi(q)} \sum_{\chi\pmod q} \biggl|\sum\Sb n\in{\Cal S}(x,y)\\
(n,q)=1\endSb \frac{\chi(n) \overline {f(n)}+1}{2} \biggr|^{2k} 
\ge \Psi(x,y;\chi_0)^{2k} e^{-3y}
\biggl(1+O\biggl(\frac{\log x}{\log q}\biggr)
\biggr)^{2k}.
$$
We deduce immediately that there are at least $\phi(q)e^{-4y} 
(1+O(\frac{\log x}{\log q}))^{2k}$ characters $\chi\pmod q$ with 
$$
\biggl|\sum\Sb n\in{\Cal S}(x,y)\\
(n,q)=1\endSb \frac{\chi(n) \overline {f(n)}+1}{2} \biggr|
\ge \Psi(x,y;\chi_0) e^{-\frac{2y}{k}} \biggl(1+O\biggl(\frac{\log x}{\log q}
\biggr)\biggr).
$$
Choosing $k=[\log q/(\log x (\log \log q)^4)]$ we conclude that 
there are $\ge q^{1-\frac{1}{(\log \log q)^2}}$ characters $\chi\pmod q$ 
for which 
$$
\biggl|\sum\Sb n\in{\Cal S}(x,y)\\
(n,q)=1\endSb \frac{\chi(n) \overline {f(n)}+1}{2} \biggr|
= \Psi(x,y;\chi_0) \biggl(1+ O\biggl(\frac{1}{(\log \log q)^4}\biggr)\biggr).
\tag{7.3}
$$
Let 
$\alpha=(\sum_{n\in{\Cal S}(x,y),\ (n,q)=1} \chi(n) \overline {f(n)})/\Psi(x,y;\chi_0)$, so that $|\alpha|\leq 1$, and (7.3) states that $|\alpha+1|=2+O(1/L^4)$ where $L=\log\log q$. By the triangle inequality
we have $2\geq 1+ |\alpha|\geq |\alpha+1|=2+O(1/L^4)$, and so 
$|1-\alpha|^2 = 2(1+|\alpha|^2)-|\alpha+1|^2=O(1/L^4)$. Thus
$|1-\alpha|=O(1/L^2)$ and the Lemma follows.
\enddemo

\demo{Proof of Theorem 3}  We suppose that $\log x \le \frac{(\log \log q)^2}{
(\log \log \log q)^2}$.
Let $y$ be as in Lemma 7.1, and put $y_1=\log q (\log \log q)^7$.  
Using Theorem 1 and Lemma 3.7, 
we get that with at most 
$q^{1-\frac{1}{\log x}}$ exceptions 
$$
\align
\sum_{n\le x}\chi(n) 
&=\Psi(x,y_1;\chi)+ O\biggl(\frac{\Psi(x,y_1)}{(\log \log q)^2}
\biggr)\\
& =\Psi(x,y;\chi) +O(|\Psi(x,y_1)-\Psi(x,y)|)
+O\biggl(\frac{\Psi(x,\log q)}{(\log \log q)^2}\biggr)\\
&= \Psi(x,y;\chi)+O\biggl( \Psi(x,\log q) \frac{\log x (\log \log \log q)^2}{(\log \log q)^2}  \biggr).\tag{7.4}\\
\endalign
$$

Given any angle $\theta$, we take $f(n)=n^{\frac{i\theta}{\log x}}$ in Lemma 
7.1.  We deduce that there are at least $q^{1-\frac{1}{\log x}}$ characters
$\chi\pmod q$ with 
$$
\align
\Psi(x,y;\chi) &= \sum\Sb n\in {\Cal S}(x,y)\\ (n,q)=1\endSb 
n^{\frac{i\theta}{\log x}} + O\biggl(\frac{\Psi(x,y;\chi_0)}{(\log \log q)^2}\biggr)\\
&= e^{i\theta} \Psi(x,y;\chi_0) + O\biggl(\sum\Sb n\in {\Cal S}(x,y) \\ (n,q)=1\endSb \frac{\log (x/n)}{\log x} + \frac{\Psi(x,y;\chi_0)}{(\log \log q)^2}\biggr)\\
&=e^{i\theta} \Psi(x,y;\chi_0) + O\biggl(\frac{\Psi(x,\log q)}{\log x}
\biggr) \\
\endalign
$$
by Corollary 3.6, and Lemma 3.7. 
Theorem 3 follows by combining this with (7.4).

\enddemo

\head 8. Results conditional on GRH: Proof of Theorem 2 
\endhead

\noindent We begin with two standard lemmas which we shall use to prove the 
conditional Theorem 2.

\proclaim{Lemma 8.1} Let $s=\sigma+it$ with $\sigma > \frac 12$ and $|t|\le 
3q$.  Let $\frac 12\le \sigma_0 < \sigma$, and 
suppose that there are no zeros of $L(z,\chi)$ inside the 
rectangle $\{ z: \sigma_0 \le \text{Re}(z) \le 1, \ \ 
|\text{Im}(z)-t|\le 3\}$.  Then
$$
|\log L(s,\chi)| \ll  \frac{\log q}{\sigma -\sigma_0}.
$$
\endproclaim
\demo{Proof} First note that if $\sigma \ge 2$ then $|\log L(s,\chi)| \ll 1$ 
and there is nothing to prove.  
We may hence assume that $\sigma <2$.  Consider the circles
with centre $2+it$ and radii $r:=2-\sigma < R:=2 -\sigma_0$, 
so that the smaller circle passes through $s$.  By our 
hypothesis $\log L(s,\chi)$
is analytic inside the larger circle.
For a point $z$ on the larger circle we use the estimate 
$|L(z,\chi)| \le 2q |z| \le q^3$,  so that 
$$
\RE \log L(z,\chi) = \log |L(z,\chi)| \le 3\log q.
$$
The Borel-Caratheodory theorem precisely states that for any point on the 
smaller circle (and so for $s$ in particular) we have
$$
\align
|\log L(s,\chi)| &\le \frac{2r}{R-r} \max_{|z-2-it|=R} \RE \log L(z,\chi) + 
\frac{R+r}{R-r} |\log L(2+it,\chi)| \\
&\ll \frac{1}{\sigma -\sigma_0}\  \log q + \frac{1}{\sigma -\sigma_0}
 \ll \frac{\log q}{\sigma -\sigma_0}.\\
\endalign
$$
\enddemo

\proclaim{Lemma 8.2}  Let $s=\sigma +it $ with $\sigma > \frac 12$
and $|t| \le 2q$.  Let $y\ge 2$ be a real number, let $\frac 12 \le 
\sigma_0 < \sigma$.  Suppose 
that there are no zeros of $L(z,\chi)$ inside the rectangle 
$\{ z: \sigma_0 \le \text{Re}(z) \le 1, \ \ |\text{Im}(z)-t|
\le y+3\}$.
Put $\sigma_1 =\min (\frac{\sigma+\sigma_0}{2}, \sigma_0+ \frac{1}{\log y})$.
Then  
$$
\log L(s,\chi) = \sum_{n=2}^{y} \frac{\Lambda(n)\chi(n)}{n^s \log n} 
+ O\biggl( \frac{\log q}{(\sigma_1 -\sigma_0)^2}
y^{\sigma_1 -\sigma}
\biggr). 
$$
\endproclaim
\demo{Proof}  Without loss of generality we may assume that 
$y\in {\Bbb Z} + \frac12$. By Perron's formula (see [3]) 
we obtain, with $c=1-\sigma+\frac 1{\log y}$,
$$
\align
\frac{1}{2\pi i} \int_{c -iy}^{c+iy } 
\log L(s+w,\chi) \frac{y^w}{w}dw 
&= \sum_{m=2}^{y}
\frac{\Lambda (m) \chi(m)}{m^s \log m} + O\biggl( \frac 1 y 
\sum_{n=1}^{\infty} \frac{y^{c}}{n^{\sigma +c}} 
\frac{1}{|\log (y/n)|} 
\biggr)\\ 
&= \sum_{m=2}^y \frac{\Lambda (m) \chi(m)}{m^s \log m} + 
O({y^{-\sigma}}\log y).
\tag{8.1}\\
\endalign
$$ 
We move the line of integration from the line Re$(w)=c$ to the 
line Re$(w)=\sigma_1 - \sigma <0$.  Our hypothesis ensures that 
the integrand is regular over the region where the line is moved, except for 
a simple pole at $w=0$ with residue $\log L(s,\chi)$.  Hence the left side
of (8.1) equals $\log L(s,\chi)$ plus
$$
\align
&\frac{1}{2\pi i} \biggl( \int_{c-iy}^{\sigma_1-\sigma-iy}
+\int_{\sigma_1-\sigma-iy}^{\sigma_1-\sigma+iy}+ 
\int_{\sigma_1-\sigma+iy}^{c +iy}\biggr) 
\log L(s+w,\chi) \frac{y^w}{w} dw \ll \frac{\log q}{(\sigma_1-\sigma_0)^2} 
y^{\sigma_1-\sigma} , \\
\endalign
$$
using Lemma 8.1 to estimate $\log L(s+w,\chi)$ in the above integrals. 
The result follows
\enddemo

If we assume the GRH for $L(s,\chi)$ then the hypotheses of Lemmas 8.1 and 
8.2 are met with $\sigma_0=\frac 12$, 
and so the conclusions drawn there are valid.  The 
advantage of these formulations is that they can be used unconditionally 
for many characters $\chi \pmod q$ by appealing to zero-density estimates; 
we exploit this to get large values of $L(\sigma,\chi)$ in [7].

We now assume the Riemann hypothesis for $L(s,\chi)$, and proceed to 
prove Theorem 2.  Define 
$$
L(s,\chi;y) = L(s,\chi) \prod_{p\le y} \biggl(1-\frac{\chi(p)}{p^s}\biggr),
$$
so that $L(s,\chi;y)$ is regular in the whole plane.  Note that
$$
\align 
\log L(s,\chi;y) &= \log L(s,\chi) + \sum_{p\le y} 
\log \biggl(1-\frac{\chi(p)}{p^s}\biggr) \\
&= \log L(s,\chi) - \sum_{m=2}^{y} \frac{\Lambda(m)\chi(m)}{m^s \log m} 
+ O\biggl(\sum\Sb p\le y\\ m\ge 2\endSb \frac{1}{mp^{m\text{Re}(s)}}\biggr),\\
\endalign
$$
and so if Re$(s)\ge \frac{1}{2} +\frac{1}{\log y}$, and 
$|\text{Im}(s)|\le 2q$
we get by Lemma 8.2
$$
|\log L(s,\chi;y)| \le C\log q \log^2 y, \tag{8.2}
$$
where $C>0$ is some constant.

Assume, 
without loss of generality, that the fractional part of $x$ is $\frac12$.
Let $u=\frac{\log x}{\log y}$ and put $c=1+\frac{1}{\log x}$.  By Perron's 
formula
$$
\align
\sum_{n\le x} \chi(n) -\Psi(x,y;\chi) &= \frac{1}{2\pi i} 
\int_{c-i\infty}^{c+i\infty} 
\biggl(L(s,\chi)-\prod_{p\le y}\biggl(1-\frac{\chi(p)}{p^s}\biggr)^{-1}
\biggr) \frac{x^s}{s} ds \\
&=\frac{1}{2\pi i} 
\int_{c-i\infty}^{c+i\infty} 
\prod_{p\le y}
\biggl(1-\frac{\chi(p)}{p^s}\biggr)^{-1} 
(\exp(L(s,\chi;y))-1) \frac {x^s}{s} ds\\ 
&=\sum_{\ell=1}^{[u]} \frac{1}{\ell!} \sum_{n\in{\Cal S}(x/y^{\ell},y)} 
\frac{\chi(n)}{2\pi i} \int_{c-i\infty}^{c+i\infty} (\log L(s,\chi;y))^{\ell} 
\L{x}{n}^s \frac{ds}{s}.\tag{8.3}\\
\endalign
$$

Now note that $(\log L(s,\chi;y))^{\ell}/\ell! = \sum_{m=1}^{\infty} 
a_{\ell}(m,y)m^{-s}$ where $|a(m,y)|\le 1$ for all $m$.  Hence, by the lemma
of section 17 of [3],
$$
\align
\frac{1}{2\pi i \ell!} \int_{c-i\infty}^{c+i\infty} 
(\log L(s,\chi;y))^{\ell} \L{x}{n}^s \frac{ds}{s} &= 
\frac{1}{2\pi i\ell!} \int_{c-ix/n}^{c+ix/n} 
(\log L(s,\chi;y))^{\ell} \L{x}{n}^s \frac{ds}{s} \\
&\hskip .75 in +O\biggl(\sum_{m=1}^{\infty}
\frac{1}{m^c} \frac{1}{|\log (x/mn)|}\biggr)  \\
&=\frac{1}{2\pi i \ell!} \int_{c-ix/n}^{c+ix/n} 
(\log L(s,\chi;y))^{\ell} \L{x}{n}^s \frac{ds}{s} +O(\log x).\\
\endalign
$$
We move the line of integration to the line segment from 
$\kappa-ix/n$ to $\kappa+ix/n$ where $\kappa:=\frac 12+\frac{1}{\log y}$.  
Using (8.2) we obtain
$$
\frac{1}{2\pi i \ell!} \int_{c-i\infty}^{c+i\infty} 
(\log L(s,\chi;y))^{\ell} \L{x}{n}^s \frac{ds}{s} 
\ll 
\L{x}{n}^{\kappa} \frac{(C\log q\log^2 y)^{\ell}}{\ell!} \log \frac{x}n 
+\log x.
$$ 
Using this in (8.3) we 
obtain
$$
\align
\biggl|\sum_{n\le x} \chi(n) -\Psi(x,y;\chi)\biggr| &\ll   
\sum_{\ell =1}^{[u]} \frac{(C\log q \log^2 y)^{\ell}}{\ell!} 
\sum_{n\in{\Cal S}(x/y^\ell,y)} \frac{x^{\kappa}}{n^{\kappa}} \log \frac{x}{n} 
+ \sum_{\ell=1}^{u} \Psi(\tfrac x{y^\ell},y) \log x.\\
\endalign
$$
Using Proposition 3.5 and Corollary 3.6 we deduce that 
(keeping in mind $y\gg \log^2 x$)
$$
\Psi(\tfrac x{y^\ell},y) \ll \biggl(\frac{c\log x}{y}\biggr)^{\ell} \Psi(x,y),
\ \ \ 
\text{and}\ \ \ 
\sum_{n\in{\Cal S}(x/y^\ell,y)} \frac{x^{\kappa}}{n^{\kappa}} \log \frac{x}{n}
\ll \ell (c\log x)^{\ell} \log^2 y\frac{\Psi(x,y)}{y^{\ell(1-\kappa)}}.
$$
Hence 
$$
\align
\biggl|\sum_{n\le x} \chi(n) -\Psi(x,y;\chi)\biggr|
&\ll \Psi(x,y) 
\sum_{\ell =1}^{[u]} \biggl(\log^2 y
\frac{(C\log q \log x \log^2 y)^{\ell}}{(\ell-1)!y^{\frac{\ell}{2}}}
+ \biggl(\frac{c\log x}{y}\biggr)^{\ell}\log x\biggr) \\
&\ll \Psi(x,y) \frac{\log q \log x \log^4 y}{y^{\frac 12}} 
\exp\biggl(O\biggl(\frac{\log q\log x \log^2 y}{y^{\frac 12}}\biggr)\biggr).
\tag{8.4}
\\
\endalign
$$
It is of interest to compare (8.4) with the bound of Theorem 6.1.

\demo{Deduction of Theorem 2}  The first assertion follows  by
taking $y=\log^2 q \log^2 x (\log \log q)^{12}$ in (8.4).  Next, 
taking $y=\log^2 q  (\log \log q)^{14}$ in (8.4) we get
$$
\sum_{n\le x} \chi(n) \ll \Psi(x,y) 
\exp\biggl(O \biggl(\frac{\log x}{(\log\log q)^3}\biggr)\biggr),
$$
and using Lemma 3.7, this is $\ll \Psi(x,\log^2 q (\log \log q)^{20})$, 
as desired. 
\enddemo

\head 9. Large character sums for real characters \endhead

\subhead 9a.  Proofs of Theorem 9, and Theorem 10 for ``small'' $x$\endsubhead

\noindent Let $y\ge 2$ be a parameter to be chosen later and put 
$b=b(y)= 4\prod_{p\le y} p$.  Choose $a\pmod b$ such that $a\equiv 1\pmod 8$,
and $\L{a}{p}=1$ for every odd $p\le y$.  Note that a squarefree integer 
$D \equiv a\pmod b$ is a fundamental discriminant satisfying $\L{D}{p}
=1$ for all $p\le y$.  We obtain the lower bounds of Theorems 9 and 10 
by averaging over fundamental discriminants of this special type, and 
choosing $y$ appropriately.

Write $n\le x$ as $n=rs$ where $p|r \implies p\le y$, and $p|s \implies 
p>y$.  Note that if $D\equiv a\pmod b$ then $\L{D}{n}=\L{D}{s}$.  Thus 
$$
\sum\Sb q\le D\le 2q \\ D\equiv a\pmod b \endSb \mu(D)^2
\sum_{n\le x} \L{D}{n} = \sum\Sb r\in{\Cal S}(x,y)\endSb \sum\Sb s\le x/r\\
p|s\implies p>y \endSb \sum\Sb q\le D\le 2q \\ D\equiv a\pmod b\endSb 
\mu(D)^2 \L{D}{s}. \tag{9.1}
$$

If $s$ is not a square then using $\mu(D)^2= \sum_{\alpha^2|D} \mu(\alpha)$
$$
\sum\Sb q\le D\le 2q\\ D\equiv a\pmod b\endSb \mu(D)^2 \L{D}{s} = 
\sum\Sb \alpha \le A \\ (\alpha, b)=1\endSb \mu(\alpha) 
\sum\Sb q\le D\le 2q\\ \alpha^2 | D\\ D\equiv a\pmod b \endSb \L{D}{s} 
+ O\biggl( \sum_{\sqrt{2q} >
\alpha > A} \biggl(\frac{q}{\alpha^2 b}+1\biggr)\biggr),
$$
and by (a modification to the proof of) the P{\' o}lya-Vinogradov inequality 
this is 
$$
\ll A \sqrt{s} \log s + \sqrt{q} + \frac{q}{A b} \ll \sqrt{q} + \frac{\sqrt{q}}
{\sqrt{b}} s^{\frac 14}\log s, \tag{9.2a}
$$
upon choosing $A= \sqrt{q}/ (b^{\frac 12} s^{\frac 14})$.  
If $s$ is a square, say $s=t^2$, then we see similarly that
$$
\sum\Sb q\le D\le 2q\\ D\equiv a \pmod b \endSb \mu(D)^2 \L{D}{s} 
= \sum\Sb q\le D\le 2q \\ D\equiv a\pmod b\\ (d,t)=1\endSb \mu(D)^2 = 
\frac{q}{b} \frac{\phi(t)}{t} \prod\Sb p > y\\ p\nmid t\endSb
 \biggl(1-\frac{1}{p^2}\biggr) 
+O(\sqrt{q} t^{\epsilon}). \tag{9.2b}
$$
Note that (9.2b) with $s=t=1$ counts the number of fundamental discriminants 
$q\le D\le 2q$ with $D\equiv a\pmod b$.

Using (9.2a,b) in (9.1) we deduce that 
$$
\!\!\!\sum\Sb q\le D\le 2q \\ D\equiv a\pmod b \endSb \!\!\mu(D)^2
\sum_{n\le x} \L{D}{n} = \frac{q}{b}
\sum\Sb r\in{\Cal S}(x,y)\endSb \sum\Sb t^2\le x/r  \\ p|t\implies p> y\endSb 
\frac{\phi(t)}{t} 
\prod\Sb p>y \\ p\nmid t\endSb \biggl(1-\frac 1{p^2}\biggr) 
+ O\biggl
(\!\sqrt{q}x^{1+\epsilon} 
+ \frac{\sqrt{q}x^{\frac 54+\epsilon}}{\sqrt{b}}\!
\biggr).\!\!\!
$$
It follows that there is at least one fundamental discriminant $D\equiv 
a\pmod b$ between $q$ and $2q$ with
$$
\sum_{n\le x} \L{D}{n} \ge \sum_{r\in{\Cal S}(x,y) } \sum\Sb t^2\le x/r  \\ p|t\implies p> y\endSb \prod_{p|t} \frac{p}{p+1} + O \biggl( \frac{b}{\sqrt q} 
x^{1+\epsilon} + \frac{\sqrt{b}}{\sqrt{q}} x^{\frac 54+\epsilon}\biggr).
\tag{9.3}
$$

We first use (9.3) to prove Theorem 9.  
Take $y$ to be the smallest prime $>\frac 13 \log q$, so that $b(y) = 
q^{\frac 13+o(1)}$.  Since 
$x\le q^{o(1)}$ we see, by counting only the $t=1$ terms on the right side of 
(9.3), that there is a fundamental discriminant $q\le D\le 2q$ with
$$
\sum_{n\le x} \L{D}{n} \ge \sum_{r\in {\Cal S}(x,y)} 1 + o(1) \ge 
\Psi(x,\tfrac 13\log q) .
$$
This proves Theorem 9.

To prove Theorem 10 in the range $\exp(\sqrt{\log q}) \leq x \leq q^{1/2}$, 
we take $y= (\frac 12 -2\epsilon) \log \frac qx$ so that
$b(y) \le (\frac qx)^{\frac 12-\epsilon}$.  By (9.3) there is a fundamental 
discriminant $q\le D\le 2q$ such that 
$$
\sum_{n\le x} \L{D}{n} \ge \sum_{r \in{\Cal S}(x,y)} 
\ \sum\Sb t^2\le x/r  \\ p|t\implies p> y\endSb \prod_{p|t} \frac{p}{p+1} 
+ O (x^{\frac 12}). \tag{9.4}
$$
We get a lower bound on the right side by counting only those $r\le R 
(\le \frac{x}{4y^2})$ for some parameter $R$ to be chosen soon.  The 
prime number theorem and the small sieve show that for such $r$ the sum 
over $t$ is $\gg \frac{\sqrt{x/r} }{\log y}$.  Hence the right side of (9.4)
is 
$$
\gg \frac{\sqrt{x}}{\log y} \sum_{r\in {\Cal S}(R,y)} \frac{1}{\sqrt{r}} 
\gg \frac{\sqrt{x}}{\log y} \frac{\Psi(R,y)}{\sqrt{R}}.
$$
Choose $R=\exp(2\sqrt {y})$ so that by Theorem 3.4 this is 
$\gg \sqrt{x} \exp( (2+o(1))\frac{\sqrt{y}}{\log y})$, as needed.

\subhead 9b. Proofs of Theorem 11, and Theorem 10 for ``large'' $x$ \endsubhead

\noindent  We 
shall consider {\sl negative } fundamental discriminants $D$, so that 
$\L{D}{-1} = -1$, and $\tau(\chi_{D}) = i\sqrt{|D|}$.  P{\' o}lya's 
Fourier expansion (see (3)) gives
$$
\align
\frac{\pi}{2\sqrt{|D|}}
\sum_{n\le |D|/N} \L{D}{n} &= \frac{1}{4} 
\sum\Sb h= -H\\ h\neq 0\endSb^{H} \frac{\L{D}{h}}{h} (1-e(-h/N)) 
+ O\biggl(\frac{1}{\sqrt{D|}}+\frac{\sqrt{|D|}}{H}\log |D|\biggr)\\ 
&= \sum_{h=1}^{H} \frac{\L{D}{h}}{h} 
\sin^2 (\pi h/N) + O\biggl(\frac{1}{\sqrt{|D|}}
+\frac {\sqrt{|D|}}{H} \log |D| \biggr).
\tag{9.5}\\
\endalign
$$

Let $y$ be a parameter to be chosen later, and let $b=b(y)$, and $a$ be
as in \S 9a.  We average (9.5) over fundamental discriminants $q \le -D \le 2q$
with $D\equiv a\pmod b$.  Arguing exactly as in the proof of (9.3), we 
deduce that there is a fundamental discriminant $D$ with $q\le -D\le 2q$ such 
that 
$$
\align
\frac{\pi}{2\sqrt{|D|} } \sum_{n\le |D|/N} \L{D}{n} 
&\ge \sum_{r\in {\Cal S}(H,y)} \frac{1}{r} \sum\Sb t^2\le H/r \\ 
p|t\implies p> y\endSb \frac{\sin^2 (\pi rt^2/N)}{t^2}
\prod_{p|t} \frac{p}{p+1} 
\\
&\hskip .5 in + O \biggl(\frac{1}{\sqrt{q}} + \frac{\sqrt{q}}{H}\log q 
+ \frac{b}{\sqrt{q}} H^{\epsilon} + \frac{\sqrt{b}}{\sqrt{q}} 
H^{\frac 14+\epsilon}\biggr).\\
\endalign
$$
Choosing $H=q^{\frac 45}/b^{\frac 25}$ we deduce that for some fundamental 
discriminant $D$ with $q\le -D\le 2q$ we have
$$
\frac{\pi}{2\sqrt{|D|} } \sum_{n\le |D|/N} \L{D}{n} 
\ge \sum_{r\in {\Cal S}(H,y)} \frac{1}{r} \sum\Sb t^2\le H/r \\ 
p|t\implies p> y\endSb \frac{\sin^2 (\pi rt^2/N)}{t^2}
\prod_{p|t} \frac{p}{p+1}  + O\biggl(q^{\epsilon}
\biggl(\frac{b}{\sqrt{q}}+\frac{b^{\frac 25}}{q^{\frac 3{10}}}\biggr)\biggr).
\tag{9.6}
$$

We now get a bound on the right side of (9.6) for various ranges of $N$.
Throughout we shall take $y=\frac 13 \log q$ so that 
$b\le q^{\frac 13+\epsilon}$.  Then $H\ge \sqrt{q}$, and the error term in 
(9.6) is $O(q^{-\frac 17})$.  

We begin with the range $\sqrt{\log q} \ge N \ge 2$.  
(Note that by taking $N=2$, the right side above is $\gtrsim 
\sum_{r\leq H} 1/r$, where the sum is over those odd $r$ whose 
prime factors are all $\leq y$, which is $\gtrsim (e^\gamma/2) 
\log y$, and we thus recover Paley's bound (4).)

We count only the terms for which $t=1$, and $r\le y$  with $\frac 14 \le 
\{ r/N \} \le \frac 34$ in (9.6).  Thus 
$$
\align
(9.6) &\ge \sum_{r\le y} \frac{\sin^2 (\pi r/N) }{r} 
+ O(q^{-\frac 17}) \ge \sum_{k=0}^{y/N} \ 
\sum_{(k +1/4)N \le r\le (k+3/4)N} 
\frac{\sin^2 (\pi r/N)}{r} + O(1) \\
&\ge \sum_{k=0}^{y/N} \frac{1}{N(k+1)} \frac 12 \biggl[ \frac N2
\biggr] +O(1) \ge \frac{1}{8} \log (y/N) + O(1) \ge \frac{1}{16}\log \log 
q+O(1).\\
\endalign
$$

Next we consider the range $\exp(\sqrt{\log q}) \ge N\ge \sqrt{\log q}$.  Here we bound (9.6) as follows:\ Let $\theta=1/\log (6\log N/\log y)$.
$$
(9.6) \ge \sum_{r \in{\Cal S}(Ny^\theta, y) } 
\frac{\sin^2 (\pi r/N)}{r} + O(q^{-\frac 17}) 
\ge \sum_{k=0}^{y^\theta} \frac{1}{N(k+1)} \ 
\sum\Sb (k +1/4)N \le n\le (k+3/4)N \\ p|n \implies p \le y \endSb \frac 12 + 
O(q^{-\frac 17}). 
$$
First we focus on the range $N < \exp((\log \log q)^2)$.
Appealing to the ``smooth numbers in short intervals estimate'', Lemma 3.8,
and Theorem 3.4 this is 
$$
\gg \sum_{k=0}^{y^\theta} \frac{1}{N(k+1)}  N 
\rho \biggl(\frac{\log (N(k+\frac 14))}{\log y}\biggr) 
\gg \theta \rho \biggl( \frac{ \log N}{\log y}+\theta\biggr) \log y,
$$
which gives the result since $\rho(u+1/\log (6u))\asymp \rho(u)$ by 
Lemma 3.9.

Next if $\exp(\sqrt{\log q}) \ge N \ge \exp((\log \log q)^2)$ 
we use Lemma 3.8, and ignore all but the $k=0$ term.  This gives 
$$
(9.6) \gg \frac 1N \Psi(\tfrac{N}{4},y) \gg \frac 1N \Psi(N,y) .
$$
The result follows from Lemma 3.9, completing the proof of Theorem 11.

To prove Theorem 10 in the range $q^{1/2}\leq x \leq q/\exp (\sqrt{\log q})$,
we consider the range $\sqrt q \ge N \ge 
\exp(\sqrt{\log q})$.  Let $R\le N/(4y^2)$ 
be a parameter to be chosen shortly.  
We bound (9.6) by considering only $r\in {\Cal S}(R,y)$, and then 
summing over values where $t=p$ is prime in the range 
$\sqrt{N}/(2\sqrt{r})\le p \le 
\sqrt{3N}/({2\sqrt{r}})$.  Thus, using the prime 
number theorem,  
$$
(9.6) \ge \sum_{r\in {\Cal S}(R,y)} \frac 1r 
\sum_{\frac{\sqrt{N}}{2\sqrt{r}
}\le p \le \frac{\sqrt{3N}}{2\sqrt{r}} } \frac{\sin^2 (\pi r/N
p^2)}{p^2} 
\gg \sum_{r \in {\Cal S}(R,y)} \frac{1}{r} \frac{\sqrt {r/N}}{\log q}
\gg \frac{1}{\sqrt{N}\log q} \frac{\Psi(R,y)}{\sqrt{R}}.
$$
Taking $R=\exp(2\sqrt{y})$ and using Theorem 3.4, this is 
$\gg (1/\sqrt{N}) \exp((2+o(1))\frac{\sqrt y}{\log y})$, as needed.

\head 10. Proof of Theorem 8 \endhead

\noindent We consider only primitive characters $\chi$ with $\chi(-1)=1$.  
Note that for a twice continuously differentiable function $\Phi$ 
the Poisson summation formula gives
$$
\sum_{n=-\infty}^{\infty} \chi(n) \Phi\L{n}{X} = \frac{X\tau(\chi)}{q}
\sum_{a=-\infty}^{\infty} \cbar(n) {\hat \Phi}\L{aX}{q}.
$$
Define $\Phi_1$ to be the characteristic function of $[-1,1]$, and 
let $\Phi_r$ be the $r$-fold convolution of $\Phi_1$.  Note that 
$\Phi_r(t)$ is supported in $[-r,r]$, $\Phi_r(-t) =\Phi_{r}(t)$, and 
that $\Phi_r(t)$ increases for $t \in [-r,0)$, and decreases for 
$t\in (0,r]$.  Lastly, note 
that ${\hat \Phi_r}(\xi) = {\hat \Phi_1}(\xi)^{r} 
= \fracwithdelims(){\sin (2\pi \xi)}{\pi \xi}^{r}$ if $\xi \neq 0$, 
and $=2^{r}$ if $\xi =0$.  We shall use the Poisson 
summation formula above with $X=q/(rN)$, and $\Phi = \Phi_r$ for an 
even value of $r\ge 4$, so that the Fourier transform ${\hat \Phi}_r$ is 
always non-negative.

On the one hand, we have 
$$
\align
\sum_{n=-\infty}^{\infty} \chi(n) \Phi_r\L{n}{X}&= 2 \sum_{n=1}^{q/N} 
\chi(n) \Phi_r\L{n}{X} = -2 \int_0^{q/N}\frac{1}{X} \Phi_{r}^{\prime}\L{t}{X} 
\sum_{n\le t} \chi (n) \ dt \\
&\le 2 \Phi_r(0) \max_{t\le q/N} \biggl| \sum_{n\le t} \chi(n)\biggr|
\le 2^r \max_{t\le q/N} \biggl|\sum_{n\le t} \chi(n) \biggr|. 
\tag{10.1}\\
\endalign
$$
On the other hand, the right side of the Poisson sum formula has size (since 
$|\tau(\chi)|=\sqrt{q}$, and $\overline{\chi}(-1)=1$)
$$
\align
\frac{2\sqrt{q}}{rN} \biggl| \sum_{a=1}^{\infty} \overline{\chi}(a) 
{\hat \Phi}_r \L{a}{rN} \biggr| 
&= \frac{2\sqrt{q}}{rN} \biggl| \sum_{a=1}^{(rN)^{\frac r{r-1}}} 
\overline{\chi}(a) \L{\sin(2\pi\frac{a}{rN})}{\frac{\pi a}{rN}}^r \biggr|
+ O\biggl( \frac{\sqrt{q}}{rN} \sum_{a> (rN)^{\frac{r}{r-1}}} 
\L{rN}{\pi a}^r  \biggr)\\
&= \frac{2\sqrt{q}}{rN} \biggl| \sum_{a=1}^{(rN)^{\frac{r}{r-1}} } 
\overline{\chi}(a) \L{ \sin(2\pi \frac{a}{rN})}{\frac{\pi a}{rN}}^r 
\biggr| + O \biggl( \frac{\sqrt{q}}{rN}\biggr). 
\tag{10.2}\\
\endalign
$$

Now observe that for integers $k \le \frac{(r-1)\log (q/2)}{r \log (rN)}$ 
we have
$$
\frac{2}{\phi(q)} \sum\Sb \chi \pmod q \\ \chi(-1)=1\endSb 
\biggl| \sum_{a=1}^{(rN)^{\frac{r}{r-1}} } 
\overline{\chi}(a) \L{ \sin(2\pi \frac{a}{rN})}{\frac{\pi a}{rN}}^r 
\biggr|^{2k} = {\Bbb E} \biggl( \biggl| \sum_{h=1}^{(rN)^{\frac{r}{r-1}}}
X_h \L{\sin (2\pi \frac{h}{rN})}{\frac{\pi h}{rN}}^{r} 
\biggr|^{2k} \biggr). 
$$
Since $r\ge 4 $ is even, note that 
$\sin^r (2\pi \frac{a}{rN})/(\frac{\pi a}{rN})^r \ge 0$ for all $a$, 
and $\ge c 2^r$ for all $a\le N$, for some absolute constant $c$.  
Hence we get from 
Lemma 2.3 that 
$$
{\Bbb E} \biggl( \biggl| \sum_{h=1}^{(rN)^{\frac{r}{r-1}}}
X_h \L{\sin (2\pi \frac{h}{rN})}{\frac{\pi h}{rN}}^{r} 
\biggr|^{2k} \biggr) \ge (c2^r)^{2k} 
{\Bbb E} \biggl( \biggl| \sum_{a=1}^{N} X_a \biggr|^{2k} \biggr).
$$
Combining the above statements thus gives
$$
\max_{t\le q/N} \biggl|\sum_{n\le t} \chi(n) \biggr| \gg
\frac{\sqrt{q}}{rN} \left( {\Bbb E} \biggl( \biggl| \sum_{a=1}^{N} X_a \biggr|^{2k} \biggr)^\frac{1}{2k} + O(1) \right) .\tag{10.3}
$$

We may obtain a lower bound from this by appealing to the results of 
\S 4 and \S 6, taking $k=[(r-1)\log (q/2)/r\log (rN)]$
in the first three parts,  choosing $r$
appropriately and replacing $x$ in those arguments by $N$ here.  
Thus the first part of the theorem is a consequence of 
Corollary 6.3 with $r=4$.  
The remaining parts of the theorem follow by choosing $r$ to 
be an even integer around $\log \log q$, and then applying
Theorem 4.1 as in the proofs of Theorems 5, 6, and 7.

\enddocument